\documentclass[11pt,a4paper]{article}

\usepackage[margin=2.7cm]{geometry}
\usepackage[T1]{fontenc}
\usepackage[utf8]{inputenc}
\usepackage{lmodern}
\usepackage{microtype}
\usepackage{amsmath,amsthm,amssymb,amsfonts,mathtools}
\usepackage[numbers,sort&compress]{natbib}
\usepackage{graphicx}
\usepackage{enumitem}%
\usepackage{aliascnt}%
\usepackage[colorlinks,linkcolor=blue,citecolor=blue,urlcolor=blue]{hyperref}
\hypersetup{pdftitle={Estimators from Convex Geometries of the Simplex},pdfauthor={Paul Navas}}

\numberwithin{equation}{section}

\theoremstyle{plain}
\newtheorem{theorem}{Theorem}[section]
\newaliascnt{lemma}{theorem}
\newtheorem{lemma}[lemma]{Lemma}
\aliascntresetthe{lemma}
\newaliascnt{proposition}{theorem}
\newtheorem{proposition}[proposition]{Proposition}
\aliascntresetthe{proposition}
\newaliascnt{axiom}{theorem}

\aliascntresetthe{axiom}
\newaliascnt{claim}{theorem}

\aliascntresetthe{claim}
\theoremstyle{definition}
\newaliascnt{definition}{theorem}
\newtheorem{definition}[definition]{Definition}
\aliascntresetthe{definition}

\newaliascnt{remark}{theorem}

\aliascntresetthe{remark}

\DeclareMathOperator{\Hess}{Hess}
\DeclareMathOperator{\diag}{diag}
\DeclareMathOperator{\clr}{clr}
\newcommand{\E}{\mathbb{E}}
\newcommand{\R}{\mathbb{R}}
\newcommand{\simplex}{\Delta^{n-1}}
\newcommand{\Dir}{\mathrm{Dir}}
\newcommand{\Var}{\mathrm{Var}}
\newcommand{\Cov}{\mathrm{Cov}}
\newcommand{\Mm}{\mathcal{M}}
\newcommand{\Ss}{\mathcal{S}}
\newcommand{\Pp}{\mathcal{P}}

\newcommand{\one}{\mathbf{1}}
\newcommand{\T}{^{\top}}

\AtBeginDocument{%
}

\begin{document}

\title{Estimators from Convex Geometries of the Simplex}
\author{Paul Navas\\[3pt]
  \normalsize Department of Statistics, TU Dortmund University\\
  \normalsize \texttt{paul.navas@tu-dortmund.de}\\
  \normalsize ORCID \href{https://orcid.org/0009-0005-6885-2108}{0009-0005-6885-2108}}
\date{\today}
\maketitle

\begin{abstract}
We unify several geometries of the probability simplex under one principle, their generation by a potential: convex potentials produce the Bregman divergences and dually flat structures of information geometry, and the exponentially concave potentials of Pal and Wong produce the logarithmic divergences and their projectively flat structures. Each structure carries a duality between maximum entropy and minimum divergence and, with it, a projection principle for estimation: the $U$-estimators of Eguchi, Komori and Ohara in the first case, the $\lambda$-deformed estimators of Wong and Zhang in the second, both governed by a Pythagorean relation and reducing to maximum likelihood in the undeformed limit. We complete the picture with estimators drawn from the third connection the divergence defines, the Levi-Civita one, by minimizing geodesic distance: existence, uniqueness and influence are then controlled by curvature and diameter rather than by flatness. Along the way the Aitchison geometry appears as the logarithmic-barrier instance of the construction at its critical weight, and the symmetrization of the Pal-Wong transport cost, a smoothing of the Hilbert projective metric, is shown to lose the dual pair and to support instead a Huber-type robust center with bounded influence in a geometry of infinite diameter.
\end{abstract}



\section{Introduction}\label{sec:Intro}

In this work, we analyze the geometry of the (interior of the) $(n - 1)$-dimensional simplex, $\simplex$, which can be thought of as the space of discrete probability distributions on $n$ atoms,

\[
p = (p_1,\dots,p_n) \in \simplex := \Bigl\{p_i > 0,\ \sum_{i=1}^np_i = 1\Bigr\}.
\]

From another perspective, seeing the set $[n]$, i.e. any set of cardinality $n$, as an underlying measurable space, with $\sigma$ algebra given by the power set $2^n$, we call $\Ss$ the set of finite measures on $([n], 2^n)$ and $\Mm$ the subset of positive measures. For any finite $m \in \Mm$ (i.e. $m([n]) < \infty$), we can take $m([n])$ as a ``num\'{e}raire'', and the ``normalized'' measure $m_n(A) \coloneqq m(A) / m([n])$ for $A \in 2^n$ lives in the simplex. If we define an equivalence relation on $\Ss$ by $s_i \sim s_j$ if $s_i(A) = k \cdot s_j(A)$, for all $A \in 2^n$, $k > 0$, i.e. the measures are equivalent up to a constant, each equivalence class is represented by a probability measure on the simplex.\\

We unify different geometric perspectives on the space of probability measures, represented on the unit simplex as above, under a general principle which generates each geometry from a convex function, its \textit{generator}\footnote{It is also known as potential, in the work of Pal and Wong, who first presented in such a generating role}. Each convex generating function furnishes a classical notion of a Bregman divergence, an instance of it being the well-known Kullback-Leibler (KL) divergence, which stems from the negative Shannon entropy generator. From the induced divergence, following a procedure recently devised by Eguchi et al in \citep{eguchikomoriohara2014}, one can derive a (Riemannian) metric and affine connections, whose geodesics deliver a statistical estimation principle that translates the geometric duality of the connections into a duality of maximum entropy and minimum divergence. A textbook description of the so-called Pythagorean theorem in this setting can be found in \citep{eguchikomori2022}.\\

A new family of divergences also stemming from convex generators was introduced by Pal and Wong in \citep{palwong2016} under the name $L$-divergences, the $L$ standing for logarithmic. There they present the well-known example of the cross-entropy being the generator of a family of divergences indexed by the points of the simplex, whereby the original convex generator is the geometric mean (whose logarithm is the said cross-entropy). In \cite{palwong2018}, again following the strategy by Eguchi in \citep{eguchikomori2022}, they construct a new Information Geometry, i.e. they approximate the $L$-divergence to second order and apply a correction term driven by the gradient of the potential, obtaining a Riemannian tensor; taking the barycenter as reference point, the tensor is the Fisher-Rao (FR) but differs from it at every other point. Following the template, from a third approximation, they obtain dual connections and show that these are \textit{projectively} dually-flat, as opposed to the dual-flatness of the Eguchi connections.\\
 
Following the template of \citep{eguchikomori2022}, a first step towards completing the ``MaxEnt-MinDiv'' program with the $L$-divergences was afforded by \citep{zhangwong2022}, who introduced a generalized Legendre duality approach, which they call a $\lambda$ deformation, based upon they derive the corresponding $\lambda$ conjugates, $\lambda$ exponential and mixture families and a conformal Hessian metric, whereby the correspondingly extended dual connections are also projectively flat, now with constant curvature $-\lambda$ (Pal and Wong's is the case $\lambda = 1$). They also connect these with Tsallis and Rényi entropy, a step which is then used by \citep{tian2026} to construct a maximum likelihood estimation method for the $\lambda$-exponential family.\\

Our main contribution completes this program with a third estimator, built on a connection not addressed in previous works, the metric or Levi-Civita connection, the average of the dual pair; it carries no dual structure, but it furnishes a natural notion of distance for the space, and we define the estimator as its Fr\'echet, i.e. distance minimizing, argument. Since that distance is defined through path length rather than by a formula, existence, uniqueness and influence of the estimator are controlled by the curvature and the diameter of the geometry rather than by its flatness, and we compare the three principles on that footing, through their influence functions. A second contribution is of a different kind: results that exist scattered in the literature are presented, as mentioned, under a common generating principle, which we then extend in two directions. One is a general Bayesian updating, obtained by mimicking the classical conjugate updating with the forward divergence, as is done with KL: it is closed under updating for every dually flat principle and only in a randomized sense for the deformed one, and its maximum-a-posteriori estimate is the divergence projection, as in Maximum-Likelihood estimation. The other is a two-parameter family of generators, the density-power potentials of the dually flat case read as generators of the deformed divergence, admissible up to an explicit frontier in $(\beta,\lambda)$ and having the density-power divergence and a tempered Kullback-Leibler divergence as its two edges. We also show that the Aitchison metric, like the classical FR one, is generated by a potential. We then analyze the symmetrization of the Pal-Wong transport cost, a uniformly bounded smoothing of the Hilbert projective metric: its induced dual structure degenerates onto the Levi-Civita connection of a rescaled Aitchison metric, so it supports neither divergence-projection principle, and the estimator it does support, its barycenter, is a robust, Huber-type modification of the compositional center. Finally, we apply the construction to difference-in-differences on categorical outcomes, where building the counterfactual composition from the three observed ones is the choice of a parallelism on the simplex: the log-ratio formula in use fixes one flat connection but leaves the metric undeclared, and different metrics report different treatment effects from the same four cells.\\

The paper is divided in four sections. In the next one, we will introduce the general notation and the two affine structures of the simplex, the Aitchison one from the perspective of Compositional Data Analysis (CoDA) and the Fisher-Rao one in its exponential coordinates, these being the main affine geometric structures upon which we will construct our estimators; the same section presents the generation of divergences, metrics and connections from convex potentials and from exponentially concave ones, the leading principle of the paper, recovers the Fisher-Rao and the Aitchison metric from their generators, and closes with the symmetrized transport cost and its degenerate geometry. Section 3 presents the estimation principles: the dually flat one of Eguchi, Komori and Ohara (EKO), the projectively flat one of Wong, Zhang, Tian and Yang (WZTY), and our main contribution, the Levi-Civita-Fr\'echet estimator, which minimizes the distance furnished by the Riemannian metric, whereby the distance minimizing curves are locally Levi-Civita geodesics; it then develops the barycenter of the symmetrized cost, the one estimator that cost does support, and compares all of them in terms of robustness and efficiency. In the final section, we present the application to difference-in-differences on categorical outcomes, where we show that the choice of geometry decides what ``parallel'' means and thus determines the calculation of the treatment effects, followed by an outlook.

\section{Affine structures of the probability simplex}\label{sec:Affine}

Two affine geometries of the simplex are in common use, and both arise from a potential. Aitchison \cite{aitchison} gave the simplex a vector space structure through the log-transform and centering, and with it a metric that now bears his name; a convex potential $\psi$ with positive definite Hessian carries the Bregman geometry with $\nabla^2\psi$ as metric, and Pal and Wong \citep{palwong2016,palwong2018}, in the same setting, show that an exponentially concave potential $\phi$ induces the metric $g = -\nabla^2\phi - \nabla\phi\,\nabla\phi\T$, restricted to the simplex tangent space. Such geometries can also be derived axiomatically, the object of information geometry since Chentsov \citep{cencov1982} (see \citep{murray} for an introduction and \citep{ay_erb} for another perspective on the two affine structures). Below, the metrics and connections are derived by the Eguchi recipe of differentiating the divergence a potential induces; the Fisher-Rao metric is recovered this way, the Aitchison geometry is generated by the logarithmic barrier at its critical weight, and the estimator of \autoref{subsec:lcf}, our main contribution, uses the Levi-Civita connection of these metrics. Calculations and proofs omitted from the text are in the Supplement; the appendix that carries each one is named in parentheses where it is used.

\subsection{Aitchison geometry of the simplex}\label{subsec:Aitch}

\subsubsection{Closure, perturbation, and powering}\label{subsec:Aitch_def}
The Aitchison distance on the simplex is defined using three operations, closure, perturbation and powering, which correspond to normalization, addition and scaling in a vector space \citep{aitchison}. Closure $C: \R_+^n \to \simplex$ turns a positive vector into a point of the simplex (a ``composition''); for $p,q \in \simplex$ and $a \in \R$,
\begin{equation}\label{eq:closure}
C(x) = \frac{x}{\sum_{k=1}^n x_k},\qquad
p \oplus q = C(p_1q_1,\dots,p_nq_n),\qquad
a \odot p = C(p_1^a,\dots,p_n^a),
\end{equation}
with $p \ominus q = p \oplus ((-1) \odot q)$; $a = 0$ gives the barycenter, the neutral element of perturbation. The \emph{centered log-ratio transform} is the map $\clr: \simplex_+ \to V_0$, where $V_0 = \{z \in \R^n: \sum_i z_i = 0\}$,
\begin{equation}\label{eq:clr}
\clr(x) = \left(\log\frac{x_1}{g(x)},\dots,\log\frac{x_n}{g(x)}\right),
\qquad g(x) = \left(\prod_{i=1}^n x_i\right)^{1/n},
\end{equation}
the geometric mean of the components of $x$. It carries $\oplus$ to $+$ and $\odot$ to scalar multiplication, so $\clr$ is a linear isomorphism onto $V_0$, the \emph{Aitchison embedding}, with inverse $z \mapsto C(e^{z_1},\dots,e^{z_n})$.

\subsubsection{Aitchison norm and distance}\label{subsec:Aitch_dist}
For $x,y \in \simplex_+$ the \emph{Aitchison inner product} is $\langle x,y\rangle_A := \langle\clr(x),\clr(y)\rangle_{\R^n}$, with norm $\|x\|_A^2 = \langle x,x\rangle_A$ and distance $d_A(x,y) := \|x \ominus y\|_A$, so that
\begin{equation}\label{eq:clr-isometrie}
d_A(x,y) = \|\clr(x) - \clr(y)\|.
\end{equation}
Thus, $\clr$ is an isometry from the Aitchison simplex onto $V_0$ with its Euclidean metric. The Aitchison geometry is therefore \emph{flat}, Euclidean: the simplex carries the linear and metric structure of $V_0$ through $\clr$, its geodesics are the $\clr$-preimages of straight lines, and the barycenter maps to the origin. In particular the geometry is unbounded: since $V_0 \cong \R^{n-1}$ is unbounded, a composition placing nearly all its mass on one part has $\clr$-image escaping to infinity, so the Aitchison diameter is infinite. This is the unbounded case, set against the finite Fisher-Rao diameter \eqref{eq:diameter} below.

\subsection{The Fisher-Rao affine structure}
We present now the exponential affine structure on the simplex, using the so-called exponential coordinate system (see \citep{amari2000}, Example~2.4 for a detailed exposition). For $p \in \Delta^{n}$ we define
\begin{equation}\label{eq:alr}
\operatorname{alr}(p) := \left(\log\frac{p_1}{p_n},\dots,\log\frac{p_{n-1}}{p_n}\right) \in \R^{n-1},
\end{equation}
a global coordinate system of the simplex.\footnote{In Compositional Data Analysis this is the additive log-ratio, or alr, coordinate system.} In these coordinates a distribution is written as
\begin{equation}\label{eq:exp-family}
p(x,\theta) = \exp\left(\sum_{k=1}^{n-1}\theta^k\delta_k(x) - \psi^*(\theta)\right),
\end{equation}
where the $\delta_k$ play the role of sufficient statistics, and
\begin{equation}\label{eq:cgf}
\psi^*(\theta) = \log\left(1 + \sum_{i=1}^{n-1}e^{\theta^i}\right) = -\log p_n
\end{equation}
ensures normalization. As in the Aitchison case one defines an addition by $p(x,\theta_1) \oplus p(x,\theta_2) = p(x,\theta_1 + \theta_2)$, which exhibits the affineness of the structure; the inverse chart is $\operatorname{alr}^{-1}(\theta) = (e^{\theta^1},\dots,e^{\theta^{n-1}},1)/Z(\theta)$ with $Z(\theta) = e^{\psi^*(\theta)}$. One is free to choose any ``num\'eraire'' for the normalization as long as it respects the multiplicative structure; centering by the geometric mean of all coordinates gives back $\clr$. The same identification is the starting point of Pistone and Shoaib \citep{pistoneshoaib2024}, whose axiomatic treatment builds the simplex as an affine space from the displacement map $s_p(q) = \log(q/p) - \E_{\mathrm{unif}}[\log(q/p)]$, whose uniform-centered case is $\clr$, and carries Aitchison's compositional calculus, Amari's dually affine geometry and the transport geometry of the simplex in a single affine language.

\subsection{Divergences and metrics from convex and exponentially concave potentials}\label{subsec:generator}
\subsubsection{Bregman divergences from convex potentials}

A divergence on $\simplex$ is a smooth function $D[p:q] \ge 0$ vanishing exactly when $p = q$; it measures discrepancy but it need be neither symmetric nor a distance. For completeness we recall the classical construction of the Bregman divergence, defined for a differentiable strictly convex potential $\psi$ by (see \citep{amari2016}, Chapter~1 for the complete derivation)
\begin{equation}\label{eq:bregman}
D_\psi[p:q] = \psi(p) - \psi(q) - \nabla\psi(q) \cdot (p - q).
\end{equation}
Taking the negative Shannon entropy $\psi(p) = \sum_i p_i\log p_i$ as potential yields $D_\psi[p:q] = \sum_i p_i\log(p_i/q_i) = D_{\mathrm{KL}}(p\|q)$, the Kullback-Leibler divergence on the simplex.

The core relationship in this setting is Legendre duality: for strictly convex $\psi$ the Legendre map sends $p$ to $\theta := \nabla\psi(p)$, and $\psi^*(\theta) = \sup_p\{\sum_{i<n}\theta^ip_i - \psi(p)\}$ is the conjugate \citep[Ch.~1]{amari2016}. For the negative Shannon entropy, with $p_n = 1 - \sum_{i<n}p_i$ substituted, the dual coordinates are the derivatives of the restricted potential, $\theta^i = \partial_{p_i}\psi - \partial_{p_n}\psi = \log(p_i/p_n) = \operatorname{alr}_i(p)$, the difference of ambient derivatives being the chain rule of the substitution, and the conjugate is the log-partition function \eqref{eq:cgf}, $\psi^*(\theta) = \log(1 + \sum_ie^{\theta^i})$: the Legendre dual coordinates are the exponential coordinates introduced above, which is what ``dual'' refers to.

\subsubsection{Logarithmic divergences from exponentially concave potentials}
We now state the family of divergences introduced by Pal and Wong \citep{palwong2016} in the context of volatility harvesting, where the generating potentials are logarithms of ``exponentially concave functions'' and their gradients are the functionally generated portfolios of stochastic portfolio theory, which solve an optimal transport problem on the simplex \citep{palwong2013,boucheywong2015,cuchieroschachermayerwong2019}. The object relevant for us is the $L$-divergence.

Given an exponentially concave function $\phi$, i.e. such that $\Phi = e^\phi$ is concave, its logarithmic $L$-divergence is \citep{palwong2018}
\begin{equation}\label{eq:l-divergence}
D_\phi[p:q] = \log\left(1 + \nabla\phi(q) \cdot (p - q)\right) - \left(\phi(p) - \phi(q)\right).
\end{equation}
Nonnegativity is exactly the content of exponential concavity: the supporting-hyperplane inequality for $\Phi = e^\phi$, divided by $\Phi(q)$ and followed by a logarithm, gives $D_\phi[p:q] \ge 0$ (\autoref{app:ldiv-calculations}). The $L$-divergence is thus the logarithmic counterpart of the Bregman divergence, with the linear term $\nabla\psi(q) \cdot (p - q)$ replaced by $\log(1 + \nabla\phi(q) \cdot (p - q))$.

As the canonical example, paralleling the role of the Kullback-Leibler divergence for the Bregman construction, take the weighted geometric mean as exponentially concave generator. Fix $\pi = (\pi_1,\dots,\pi_n) \in \simplex$ and set $\Phi(p) = \prod_ip_i^{\pi_i}$, so that $\phi(p) = \sum_i\pi_i\log p_i$ and $\partial_i\phi(p) = \pi_i/p_i$. Substitution into \eqref{eq:l-divergence} yields
\begin{equation}\label{eq:ldiv-amgm}
D_\phi[p:q]
= \log\left(\sum_{i=1}^n\pi_i\frac{p_i}{q_i}\right) - \sum_{i=1}^n\pi_i\log\frac{p_i}{q_i}
= \log\operatorname{AM}_\pi(r) - \log\operatorname{GM}_\pi(r),
\qquad r_i := \frac{p_i}{q_i}.
\end{equation}
This is the logarithm of the ratio of the $\pi$-weighted arithmetic mean to the $\pi$-weighted geometric mean of the likelihood ratios. Nonnegativity is here exactly the arithmetic-geometric mean inequality, with equality if and only if all $r_i$ coincide, i.e. $p = q$, recovering the general property of \eqref{eq:l-divergence} (\autoref{app:ldiv-calculations}).

\subsubsection{The metric and connections induced by a convex potential}\label{subsubsec:metric-convex}

To obtain a Riemannian metric and a pair of affine connections from a divergence, we follow Eguchi \citep[Thm.~2.6]{eguchikomori2022}, going back to \citep{eguchi1992}, who gets them by differentiating the divergence in its two arguments and evaluating on the diagonal. The statistical prototype is the expansion $D_{\mathrm{KL}}(p_\theta\|p_{\theta+d\theta}) = \tfrac{1}{2}\,d\theta\T I(\theta)\,d\theta + O(|d\theta|^3)$, Fisher information as the quadratic term of the Kullback-Leibler divergence; the recipe is this expansion for an arbitrary divergence, the metric being the second-order term and the two connections the two ways of reading the third-order term, according to which argument carries two of the three derivatives. Throughout, $D[\cdot:\cdot]: \simplex \times \simplex \to \R_{\ge0}$ is a divergence and we write $\partial_i = \partial/\partial p_i$ and $\partial_{j'} = \partial/\partial q_j$. The metric and two connections are
\begin{align}
g_{ij}(p)&:= -\left.\partial_i\partial_{j'}D[p:q]\right|_{q=p}, \label{eq:eguchi-metric}\\
\Gamma^{(m)}_{ij,k}(p)&:= -\left.\partial_i\partial_j\partial_{k'}D[p:q]\right|_{q=p}, \label{eq:eguchi-m-connection}\\
\Gamma^{(e)}_{ij,k}(p)&:= -\left.\partial_{i'}\partial_{j'}\partial_kD[p:q]\right|_{q=p}, \label{eq:eguchi-e-connection}
\end{align}
where $T_p\simplex = \{v \in \R^n: \sum_i v_i = 0\}$ and $\Gamma_{ij,k} = g_{k\ell}\Gamma^\ell_{ij}$. The names follow the statistical convention of Amari \citep{amari2000}. For the \textit{categorical family}, the saturated model of all distributions on $n$ atoms whose parameter space is the simplex itself, the mixture connection is flat in the expectation coordinates $p$, while the exponential connection is flat in the natural coordinates $\theta = \operatorname{alr}(p)$. These two connections are \emph{dual with respect to $g$}, $Xg(Y,Z) = g(\nabla^{(e)}_X Y,Z) + g(Y,\nabla^{(m)}_X Z)$, and both are torsion-free, so their average is the unique torsion-free $g$-compatible connection, i.e. the Levi-Civita connection \citep[Thm.~2.7]{eguchikomori2022},
\begin{equation}\label{eq:lc-average}
\nabla^{(\mathrm{LC})} = \tfrac{1}{2}\left(\nabla^{(e)} + \nabla^{(m)}\right);
\end{equation}
it is this metric connection whose geodesics define the distance underlying our Levi-Civita-Fr\'echet (LCF) estimator below.

We now evaluate these objects for the Bregman divergence \eqref{eq:bregman}. Differentiating once in each argument gives $\partial_iD_\psi[p:q] = \partial_i\psi(p) - \partial_i\psi(q)$ and hence
\begin{equation}\label{eq:hesse-metrik}
g_{ij}(p) = \partial_i\partial_j\psi(p),\qquad g(p) = \nabla^2\psi(p).
\end{equation}
The metric induced by a convex potential is therefore exactly its Hessian: convexity makes it positive semidefinite, strict convexity on $T_p\simplex$ a genuine Riemannian metric. Moreover, $\partial_i\partial_jD_\psi[p:q] = \partial_i\partial_j\psi(p)$ no longer depends on $q$, so $\Gamma^{(m)}_{ij,k}(p) = 0$: the potential's coordinates $p$ are affine for $\nabla^{(m)}$. By duality, $\nabla^{(e)}$ is flat in the Legendre-conjugate coordinates $\theta = \nabla\psi(p) = \operatorname{alr}(p)$. The Bregman geometry is therefore \emph{dually flat}, the property on which the EKO estimation principle of \autoref{sec:Estimators} rests.

We specialize to the same potential that produced the Kullback-Leibler divergence, the negative Shannon entropy, with Hessian $\partial_i\partial_j\psi(p) = \delta_{ij}/p_i$. By \eqref{eq:hesse-metrik},
\begin{equation}\label{eq:fr-metrik}
g_{ij}(p) = \frac{\delta_{ij}}{p_i},\qquad
g_p(v,w) = \sum_{i=1}^n\frac{v_iw_i}{p_i},\quad v,w \in T_p\simplex,
\end{equation}
the Fisher-Rao metric of the categorical family. Its mixture connection vanishes in $p$ and its exponential connection in $\theta = \operatorname{alr}(p)$; the Levi-Civita connection \eqref{eq:lc-average} is their average and is, in general, nonzero. Its curvature is exhibited by the square-root embedding
\begin{equation}\label{eq:sphere-iso}
\Pi: \simplex \longrightarrow S_2^{n-1} \cap \R^n_{>0},\qquad
\Pi(p) = 2(\sqrt{p_1},\dots,\sqrt{p_n}),
\end{equation}
onto the open positive orthant of the sphere of radius $R = 2$. Writing $\xi_i = 2\sqrt{p_i}$ gives $d\xi_i = p_i^{-1/2}dp_i$, so the round metric pulls back to $\sum_i dp_i^2/p_i$, precisely \eqref{eq:fr-metrik}. Hence $\Pi$ is an isometry onto its image and the Fisher-Rao metric has constant sectional curvature $\kappa = 1/R^2 = 1/4$. The dually flat affine structure and the curved Levi-Civita structure therefore coexist on the same simplex: the connection one selects, and not the metric alone, determines the geodesics.

The \emph{geodesic distance} $d(p,q) := \inf_\gamma\int_0^1\sqrt{g_{\gamma(t)}(\dot\gamma(t),\dot\gamma(t))}\,dt$, the infimum over curves joining $p$ and $q$, is the quantity whose square is the LCF loss; a minimizing curve is a Levi-Civita geodesic, and under \eqref{eq:sphere-iso} the distance is the great-circle arc length between $\Pi(p)$ and $\Pi(q)$.

The \emph{geodesic diameter}, the largest distance the geometry allows,
\begin{equation}\label{eq:diameter}
\operatorname{diam}(\simplex,d) := \sup_{p,q\in\simplex}d(p,q) \le R\frac{\pi}{2} = \pi,
\end{equation}
because two points in one orthant meet at an angle of at most $\pi/2$. Thus Fisher-Rao has finite diameter: no two distributions in the family are more than a bounded distance apart, whereas Aitchison, isometric to Euclidean space, has infinite diameter. Finite diameter is the geometric source of the bounded influence of the Fisher-Rao Fr\'echet estimator, made precise in \autoref{subsec:robustness}.

\subsubsection{The metric and connections induced by an exponentially concave potential}\label{subsubsec:metric-lconvex}

The Eguchi construction \eqref{eq:eguchi-metric}--\eqref{eq:lc-average} takes \emph{any} divergence as input, so we may apply it directly to the logarithmic divergence \eqref{eq:l-divergence}.

Writing $A(p,q) := 1 + \nabla\phi(q) \cdot (p - q)$ and $g^{\mathrm H} := -\nabla^2\phi$ for the ``bare-Hessian'' part, differentiation of \eqref{eq:l-divergence} in its two arguments gives, on the diagonal (\autoref{app:ldiv-eguchi}),
\begin{equation}\label{eq:ldiv-metric}
g = -\left(\nabla^2\phi + \nabla\phi\,\nabla\phi\T\right)
= g^{\mathrm H} - \nabla\phi\,\nabla\phi\T
= -\Phi^{-1}\nabla^2\Phi,\qquad \Phi = e^\phi,
\end{equation}
the last equality by the chain rule $\nabla^2\Phi = \Phi(\nabla^2\phi + \nabla\phi\nabla\phi\T)$. Concavity of $\Phi$ makes $g$ positive semidefinite: positive semidefiniteness is exponential concavity restated. The metric is the Hessian metric corrected by a rank-one term built from the gradient of the potential, the geometric contribution of the logarithm in \eqref{eq:l-divergence}, absent in the Bregman case, where \eqref{eq:hesse-metrik} gave the Hessian alone.

The same one-order-higher differentiation produces the connections \eqref{eq:eguchi-e-connection}--\eqref{eq:eguchi-m-connection}, again with Levi-Civita average \eqref{eq:lc-average}. Unlike its Bregman counterpart, the relevant derivative still depends on both arguments, so $\nabla^{(m)}$ is not flat in the coordinates $p$ (\autoref{app:ldiv-eguchi}); what survives is projective flatness, the central result of Pal and Wong \citep{palwong2018}, stated as \autoref{prop:proj-flat} where it drives the WZTY estimation principle.

As example for this exponentially concave setting, take the weighted geometric mean $\Phi(p) = \prod_i p_i^{\pi_i}$ with fixed weights $\pi \in \simplex$, whose logarithm $\phi(p) = \sum_i\pi_i\log p_i$ is the negative cross-entropy. With $\partial_i\phi = \pi_i/p_i$ and $\partial_i\partial_j\phi = -\delta_{ij}\pi_i/p_i^2$, \eqref{eq:ldiv-metric} gives
\begin{equation}\label{eq:ldiv-beispiel-metrik}
g_{ij}(p) = \frac{\pi_i}{p_i^2}\delta_{ij} - \frac{\pi_i\pi_j}{p_ip_j}.
\end{equation}
The diagonal part is the bare Hessian; the rank-one part is the correction \eqref{eq:ldiv-metric}. To identify its Riemannian geometry, define the weighted analogue of the centered log-ratio isometry,
\begin{equation}\label{eq:clr-pi}
\clr_\pi(p)_i = \log p_i - \sum_{j=1}^n\pi_j\log p_j,
\qquad
V_\pi = \left\{z \in \R^n: \sum_i\pi_i z_i = 0\right\},
\end{equation}
with constant inner product $\langle a,b\rangle_\pi = \sum_i\pi_i a_i b_i$. The weighted-variance calculation gives (\autoref{app:weighted-clr})
\begin{equation}\label{eq:ldiv-clr-rueckzug}
g = \clr_\pi^*\langle\cdot,\cdot\rangle_\pi.
\end{equation}
A pullback of a flat metric under a diffeomorphism is flat, so the $L$-divergence metric of the weighted geometric mean has vanishing Riemannian curvature for every weight $\pi$. For uniform weights, $\clr_\pi = \clr$ and $\langle\cdot,\cdot\rangle_\pi = \tfrac{1}{n}\langle\cdot,\cdot\rangle_{\R^n}$, hence $g = \tfrac{1}{n} g^A$.

At the barycenter with uniform weights, \eqref{eq:ldiv-beispiel-metrik} gives $g_{ij} = n\delta_{ij} - 1$, so $g_p(v,v) = n\sum_iv_i^2$ and $g_p = g^{FR}_p$ as bilinear forms; the two geometries built from the same logarithmic potentials agree at that one point and nowhere else, one curving onto the sphere of radius $2$, the other flat.

\subsubsection{The generators for the Aitchison geometry}\label{subsubsec:aitchison-generators}

The two construction principles, Bregman divergences from convex potentials and $L$-divergences from exponentially concave ones, intersect naturally on a single family of potentials, the logarithmic generators, the Aitchison geometry of \autoref{subsec:Aitch_dist} appearing as an instance rather than as a separate construction. Fix $c > 0$ and consider the logarithmic barrier
\begin{equation}\label{eq:logbarriere}
\psi_c(p) = -c\sum_{i=1}^n\log p_i,\qquad
\nabla^2\psi_c(p) = c\,\diag(p_1^{-2},\dots,p_n^{-2}).
\end{equation}
It is strictly convex, so the convex recipe \eqref{eq:hesse-metrik} returns the bare-Hessian metric
\begin{equation}\label{eq:barrier-metric}
g^{\mathrm{bare}}_c{}_p(u,u) = c\sum_{i=1}^n\frac{u_i^2}{p_i^2},\qquad u \in T_p\simplex.
\end{equation}
The denominator is $p_i^2$, one power higher than the $p_i$ of Fisher-Rao: the barrier supplies the Aitchison \emph{scale}, measuring relative log-ratio change, but not yet the Aitchison geometry, because \eqref{eq:barrier-metric} is curved.

The log chart $y_i = \log p_i$ carries \eqref{eq:barrier-metric} to the metric induced by $c\sum_idy_i^2$ on the hypersurface $\sum_ie^{y_i} = 1$. The defining function has positive definite Hessian, so the hypersurface is strictly convex, its second fundamental form is definite, and by the Gauss equation every sectional curvature of the induced metric is positive for $n \ge 3$, a sign the rescaling by $c$ preserves. The barrier does not produce Aitchison as a convex potential; it does so as an exponentially concave one.

Taking $\phi(p) = c\sum_i\log p_i$, exponential concavity of $e^\phi = (\prod_i p_i)^c$ holds precisely for $c \le 1/n$, and \eqref{eq:ldiv-metric} gives
\begin{equation}\label{eq:barrier-Ldiv-metric}
g^L_c{}_p(u,u) = c\sum_{i=1}^n\frac{u_i^2}{p_i^2} - c^2\left(\sum_{i=1}^n\frac{u_i}{p_i}\right)^2
= c\left[nc\,g^A_p(u,u) + (1 - nc)g^{\mathrm{bare}}_1{}_p(u,u)\right].
\end{equation}
This is a $c$-scaled interpolation between the flat Aitchison metric and the curved barrier metric. The barrier component carries the weight $1 - nc$ and disappears exactly at the upper exponential-concavity threshold $c = \tfrac{1}{n}$; the rank-one correction, which distinguishes the exponentially concave perspective from the convex one, is what removes the curvature. There the metric collapses to $g^L_{1/n} = \tfrac{1}{n} g^A$, the uniform-weight case of \eqref{eq:ldiv-clr-rueckzug}. For $c < 1/n$ the barrier component survives; flatness occurs at $c = 1/n$, where exponential concavity becomes tight, and, for $n \ge 3$, nowhere else.\footnote{The exact Aitchison metric can be recovered by deforming, rather than rescaling, the unscaled barrier at $\lambda = 1/n$; see \autoref{subsec:wzty}.}

\subsubsection{The symmetrized cost and its degenerate geometry}\label{subsubsec:symcost}

The $L$-divergence, like every divergence, has two orientations, and $D_\phi[p:q] \ne D_\phi[q:p]$. Their sum\footnote{The $L$-analogue of passing from Kullback-Leibler to Jeffreys' divergence.} carries a second reading: for exponentially concave $\phi$ the divergence is the cost function of an optimal transport problem on the simplex, whose optimal maps are generated by the potential \citep{palwong2018,palwong2020}, so symmetrizing a divergence is here also symmetrizing a transport cost. We do so for the uniform-weight generator $\phi(p) = \frac{1}{n}\sum_i\log p_i$ of Pal and Wong, for which the two orientations carry the same log-ratio term with opposite signs, so the symmetrization retains only the two arithmetic means of the likelihood ratios:
\begin{equation}\label{eq:symkosten}
D^{\mathrm{sym}}[p:q] := D_\phi[p:q] + D_\phi[q:p]
= \log\left(\frac{1}{n}\sum_i\frac{p_i}{q_i}\right) + \log\left(\frac{1}{n}\sum_i\frac{q_i}{p_i}\right).
\end{equation}
This is a symmetric divergence, nonnegative by Cauchy-Schwarz, with equality only at $p = q$. In the centered log-ratio chart the cost is a function of the difference alone: writing $\Lambda(v) = \log(\frac{1}{n}\sum_i e^{v_i})$, $\xi = \clr(p)$ and $\eta = \clr(q)$, the geometric-mean factors cancel between the two terms of \eqref{eq:symkosten}, leaving
\begin{equation}\label{eq:symkosten-clr}
D^{\mathrm{sym}}[p:q] = h(\xi - \eta),\qquad h(v) = \Lambda(v) + \Lambda(-v),\quad v \in V_0,
\end{equation}
an even function of the clr difference. Three structural properties follow, each needed below (proofs in \autoref{app:symmetry}).

First, invariance. $D^{\mathrm{sym}}$ is invariant under perturbation, permutations of the parts, the swap of its arguments and the \emph{compositional inversion}
\begin{equation}\label{eq:inversion}
\iota(p) = C(p_1^{-1},\dots,p_n^{-1}),\qquad \clr(\iota(p)) = -\clr(p),
\end{equation}
an involution fixing the barycenter. The one-sided divergence is not $\iota$-invariant; the inversion exchanges its orientations, $D_\phi[\iota(p):\iota(q)] = D_\phi[q:p]$ for the uniform-weight generator, which is what makes the sum invariant.

Second, convexity. The Hessian of $h$ is
\begin{equation}\label{eq:symcost-hessian}
\nabla^2h(v) = \Cov_{p} + \Cov_{\iota(p)},
\qquad p = \clr^{-1}(v),
\qquad \Cov_q = \diag(q) - qq\T,
\end{equation}
the covariance matrices of the categorical law $p$ and of its inversion \eqref{eq:inversion}, both strictly positive, so both summands are positive semidefinite with kernel exactly $\operatorname{span}(\one)$: $h$ is strictly convex on $V_0$, at every point, not only near the origin. The same matrices written through the exponential weights of $\pm v$ are in \autoref{app:hessian-weights}. Third, growth. To assess it, we use the Hilbert projective metric $d_{\mathrm{Hil}}(p,q) = \log\max_i(p_i/q_i) - \log\min_i(p_i/q_i)$, which in clr coordinates is the range $\max_i v_i - \min_i v_i$ of the difference.

\begin{lemma}[Hilbert lower and upper bounds]\label{lem:hilbert-sandwich}
For all $p,q \in \simplex$,
\begin{equation}\label{eq:hilbert-sandwich}
\max(0,d_{\mathrm{Hil}}(p,q) - 2\log n) \le D^{\mathrm{sym}}[p:q]
\le 2\log\cosh\left(\frac{d_{\mathrm{Hil}}(p,q)}{2}\right) \le d_{\mathrm{Hil}}(p,q).
\end{equation}
\end{lemma}

The symmetrized cost is therefore a smooth approximation to the Hilbert metric: it stays within $2\log n$ of it everywhere, but replaces the corners and flat linear regions of the Hilbert range function by a smooth, strictly convex profile. Near the diagonal it is quadratic: since $h$ is even, all odd-order terms in its Taylor expansion vanish, and
\begin{equation}\label{eq:symcost-germ}
h(v) = \frac{1}{n}|v|^2 + O(|v|^4),\qquad v \in V_0.
\end{equation}
Thus the cost grows quadratically near the diagonal and linearly, at the Hilbert slope, in the tails, i.e. the profile of a Huber loss, with the crossover between the two regimes fixed by the cost itself rather than chosen by hand. The statistical use of this shape is the subject of \autoref{subsec:jbary}. The geometry the Eguchi construction extracts from the divergence is degenerate: for any symmetric divergence the two Eguchi connections coincide and reduce to the Levi-Civita connection \citep[p.~41]{eguchikomori2022}, since the connections are third derivatives at the diagonal and odd under the swap of arguments, so symmetrization kills the cubic tensor $\Gamma^{(e)} - \Gamma^{(m)}$.

\begin{proposition}[The Eguchi geometry of the symmetrized cost]\label{prop:symcost-degenerate}
The metric induced by $D^{\mathrm{sym}}$ through \eqref{eq:eguchi-metric} is
\begin{equation}\label{eq:symcost-metric}
g^{\mathrm{sym}} = \frac{2}{n} g^A,
\end{equation}
twice the uniform-weight $L$-divergence metric of \eqref{eq:ldiv-beispiel-metrik}, hence flat, and in the clr chart both connections \eqref{eq:eguchi-m-connection}--\eqref{eq:eguchi-e-connection} vanish identically, so the primal and dual flat charts coincide.
\end{proposition}

The cancellation is concrete here: each orientation of the $L$-divergence contributes a projective term to its mixture connection and its dual contributes the opposite one; adding the two orientations cancels them. What remains is the Aitchison structure, at the critical weight of \autoref{subsubsec:aitchison-generators}, up to the factor $2$. The Eguchi construction builds its metric from second derivatives at the diagonal and its connections from the third, so it sees only how $D^{\mathrm{sym}}$ behaves between nearby points. The bounded distortion from Hilbert, the linear tails and the Huber profile are global features and cannot be seen by this construction.

One could still hope that the divergence belongs somehow to one of the two classes on which \autoref{sec:Estimators} rests, convex or exponentially concave, under some other chart or generator. It does not, as the next proposition shows.

\begin{proposition}[Exclusion from both divergence classes]\label{prop:symcost-exclusion}
$D^{\mathrm{sym}}$ is not the Bregman divergence of any strictly convex potential in any smooth chart of the simplex, and not the $\lambda$-logarithmic divergence \eqref{eq:lambda-ldiv} of any $\lambda$-exponentially concave generator for any $\lambda > 0$.
\end{proposition}

Neither projection principle therefore applies to $D^{\mathrm{sym}}$. We instead minimize it directly as a loss, a Fr\'echet-type principle applied to the divergence rather than to the squared distance of its metric, and the two differ: the estimator is the barycenter of \autoref{subsec:jbary}, while the Levi-Civita-Fr\'echet mean of the induced metric is the geometric mean of case (b) in \autoref{subsec:lcf}.

\section{General estimation principles}\label{sec:Estimators}

\subsection{The estimation problem}\label{subsec:estimation}

With the mathematical tools at hand, we now see the simplex also as sample space. A point of $\simplex$ is a probability assignment over $n$ alternatives, which can be either observed categories (market shares or voting proportions) or the elements of a finite parameter set $\Theta = \{\theta_1,\dots,\theta_n\}$\footnote{in the latter case it is a belief over $\Theta$, a prior, a posterior, or the normalized likelihood vector $\ell = C(f(y \mid \theta_1),\dots,f(y \mid \theta_n))$ of one observation.}. In either case we observe independent $x_1,\dots,x_N \sim \Pp$ on $\simplex$ and seek a representative selected by a loss supplied by the geometry. Maximum likelihood and maximum entropy will reappear below as the classical instances of the EKO projection principle (\autoref{subsubsec:likelihood} and \autoref{subsec:eko}), while the CoDA center appears as the flat Levi-Civita-Fr\'echet mean in \autoref{subsec:lcf}.

We use repeatedly the arithmetic mean $\bar x = \frac{1}{N}\sum_{k=1}^N x_k$. A categorical sample, encoded at the vertices by $x_k = e_{X_k}$, has $\bar x = \hat p$, the empirical distribution $\hat p_i = \#\{k: X_k = i\}/N$. The theory is developed for interior observations, the categorical case following wherever the loss extends finitely to the vertices; the one geometry for which it does not is taken up in \autoref{subsec:lcf}.

We restrict candidate representatives to a \emph{model}.
\begin{definition}[Model]\label{def:model}
A \emph{model} is a subset $M \subseteq \simplex$ of candidate points in the manifold. For our estimators we use the moment-constraint family
\begin{equation}\label{eq:model}
  M_\tau = \{p \in \simplex: \E_p[t] = \tau\},
  \qquad
  t: \{1,\dots,n\} \to \R^k,
  \qquad
  \E_p[t] = \sum_{i=1}^n p_i\,t(i),
\end{equation}
which is affine in the mixture coordinates. The complementary $U$-model introduced below is affine in the corresponding exponential coordinates; their intersection is the setting of the Pythagorean projection of \autoref{thm:pyth}.
\end{definition}

The moment surfaces enter below in three roles: as the mixture-affine side of the Pythagorean projection, as the constraint set of the maximum-entropy problem, and as a model over which the Fr\'echet estimator can be minimized. In general an estimator assigns to the sample a member of the model minimizing an empirical loss,
\begin{equation}\label{eq:m-estimator}
  T(x_1,\dots,x_N) = \arg\min_{q\in M}\frac{1}{N}\sum_{k=1}^N\ell(q,x_k),
\end{equation}
for a loss $\ell: \simplex \times \simplex \to [0,\infty]$; this is the empirical counterpart of the risk $\E_\Pp[\ell(q,X)]$. In decision-theoretic language the minimized risk is the generalized entropy of the loss and the excess its divergence, the correspondence underlying the maximum-entropy/minimum-divergence duality below \citep{grunwalddawid2004}, see also \citep[Thm.~3.5]{eguchikomori2022}.

For categorical data the empirical risk depends on the sample only through the empirical distribution, $\frac{1}{N}\sum_{k=1}^N\ell(q,e_{X_k}) = \sum_i\hat p_i\,\ell(q,e_i) $ (for interior observations this need not hold). Bregman losses are distinguished by an analogous reduction to the arithmetic mean, which we state in weighted form because the same identity yields the conjugacy of \autoref{prop:u-conjugacy}.

\begin{lemma}[Two-stage identity]\label{lem:two-stage}
Let $D_\psi$ be the Bregman divergence \eqref{eq:bregman} of a differentiable convex $\psi$. For any $x_1,\dots,x_N \in \simplex$, weights $w_1,\dots,w_N \ge 0$ with $W = \sum_{k=1}^Nw_k > 0$, and $q \in \simplex$,
\begin{equation}\label{eq:two-stage}
  \sum_{k=1}^N w_k D_\psi[x_k:q]
  =
  \sum_{k=1}^N w_k D_\psi[x_k:\bar x_w]
  +
  W D_\psi[\bar x_w:q],
  \qquad
  \bar x_w = \frac{1}{W}\sum_{k=1}^N w_kx_k.
\end{equation}
\end{lemma}

\begin{proof}
Averaging $D_\psi[x:q] = \psi(x) - \psi(q) - \nabla\psi(q) \cdot (x - q)$ gives
$\sum_k w_kD_\psi[x_k:q] = \sum_k w_k\psi(x_k) - W\psi(q) - W\nabla\psi(q) \cdot (\bar x_w - q)$.
At $q = \bar x_w$ the gradient term vanishes; subtracting the resulting identity gives \eqref{eq:two-stage}.
\end{proof}

Thus, for uniform weights, $\arg\min_{q\in M}\frac{1}{N}\sum_{k=1}^ND_\psi[x_k:q] = \arg\min_{q\in M}D_\psi[\bar x:q]$, so $\bar x$ is sufficient for the Bregman projection problem. For categorical observations $x_k = e_{X_k}$, it is precisely the empirical distribution, $\bar x = \hat p$.\footnote{For fixed $N$, $\hat p$ is equivalent to the count vector, the usual sufficient statistic for the categorical parameter. More generally, the mean-optimality property characterizes Bregman losses under mild regularity conditions \citep{banerjeeguowang2005}.} The mean is also sufficient for the $U$-divergences below, of which Bregman is an instance \citep[Thm.~4.3]{eguchikomori2022}. The $L$-divergence and the geodesic loss below do not share this, and their estimators use more information than just the mean.

The three principles of this section are three choices of $\ell$ in \eqref{eq:m-estimator}. EKO takes a $U$-divergence, for which \autoref{lem:two-stage} reduces the empirical risk to a divergence from $\bar x$; over the exponential family $E_t$ and with the Kullback-Leibler loss the projection $\arg\min_{q\in E_t}D_\psi[\bar x:q]$ is maximum likelihood, and the fit is the intersection of $E_t$ with the moment surface $M_{\bar t}$ (\autoref{subsubsec:likelihood}, \autoref{subsec:eko}). WZTY takes the $\lambda$-logarithmic divergence on the $\lambda$-exponential family $E_t^\lambda$ and keeps the projection template but deforms the duality: the sample is projected, not its mean, the score becomes an estimate-dependent reweighting (\autoref{subsec:wzty}). The third principle, our main contribution, takes the squared geodesic distance of the Levi-Civita connection as loss, yielding the Levi-Civita-Fr\'echet estimator of \autoref{subsec:lcf}; since $d^2$ is defined through path length rather than by an algebraic formula, existence, uniqueness and robustness depend on the metric geometry itself (\autoref{subsec:robustness}). For the symmetrized cost the dual pair collapses onto one connection and there is no projection principle; its direct barycenter is treated in \autoref{subsec:jbary}.

\subsubsection{Likelihood, counts, and Gibbs conditioning}\label{subsubsec:likelihood}

The classical logarithmic loss is the benchmark for the geometric losses below. For categorical data, $\ell(q,e_i) = -\log q_i$ and $-\sum_i\hat p_i\log q_i = D_{\mathrm{KL}}(\hat p\|q) + H(\hat p)$, so maximizing the likelihood of the counts over a model $M$ is the forward Kullback-Leibler projection \citep[eq.~(5)]{eguchi2021}
\begin{equation}\label{eq:ml-projektion}
  \hat q_{\mathrm{ML}}
  = \arg\max_{q\in M}\prod_{i=1}^nq_i^{N\hat p_i}
  = \arg\min_{q\in M}D_{\mathrm{KL}}(\hat p\|q):
\end{equation}
maximum likelihood is already the Shannon instance of the EKO construction.\footnote{The forward divergence also has a probabilistic reading. Sanov's theorem gives $D_{\mathrm{KL}}(q\|P)$ as the rate function for the empirical distribution of an i.i.d.\ sample from $P$ \citep{dembozeitouni1998}, and Gibbs conditioning says that, conditional on the atypical event $\hat p_N \in M_\tau$, the law of one observation converges to the $I$-projection of $P$ onto $M_\tau$ \citep{csiszar1984,vancampenhoutcover1981}, the maximum-entropy member of the surface when $P$ is uniform: the point that \autoref{subsec:eko} obtains from Pythagorean projection and moment matching.}

\subsection{Dually flat estimation: the EKO principle}\label{subsec:eko}

In the Shannon case $\psi(p) = \sum_i p_i\log p_i$ gives $\partial_i\psi(p) = \log p_i + 1$: the dual coordinate is $\theta_i = \log p_i$, and $p_i = e^{\theta_i}$ recovers the probability from it. The principle of Eguchi, Komori and Ohara (EKO) \citep{eguchikomoriohara2014} replaces this pair by $(\xi,u)$ for a scalar strictly convex increasing $U: \R \to \R$ with $u = U'$ and $\xi = u^{-1}$: $\xi$ takes the role of $\log$ and $u$ that of $\exp$. They reproduce around $(\xi,u)$ the dually flat structure of divergence, entropy, exponential-type family and Pythagorean projection; $U = \exp$ recovers Shannon, with $(\xi,u) = (\log,\exp)$. We use their construction in this general form, in the notation of the monograph \citep{eguchikomori2022}, whose Chapters 2--4 we cite for the statements we take over.

\begin{align}
  D_U[p:q]&= C_U(p,q) - H_U(p)
  = \sum_{i=1}^n\Big[U(\xi(q_i)) - U(\xi(p_i)) - p_i\big(\xi(q_i) - \xi(p_i)\big)\Big]\notag\\
  &= D_{\psi_U}[p:q],
  \qquad
  \psi_U(p) = \sum_{i=1}^nU^*(p_i),\label{eq:u-divergenz}
\end{align}
with $U^*$ the scalar Legendre conjugate of $U$. Each summand is a scalar Bregman divergence of $U$, so nonnegativity follows termwise from convexity, with equality only at $p = q$; at $U = \exp$ the affine terms cancel by normalization and $D_U$ becomes $D_{\mathrm{KL}}$, while $H_U$ is Shannon entropy up to an additive constant. The last equality says that the construction is a Bregman divergence on the simplex itself, with potential $\psi_U = -H_U$, since $(U^*)' = \xi$ gives $\nabla\psi_U(p) = (\xi(p_i))_i$ and $\nabla^2\psi_U(p) = \diag(\xi'(p_i)) \succ 0$. Hence the Hessian recipe of \autoref{subsubsec:metric-convex} gives
\begin{equation}\label{eq:u-metrik}
  g^U_{ij}(p) = \xi'(p_i)\delta_{ij},
\end{equation}
the Fisher-Rao metric at $U = \exp$; the mixture connection is flat in $p$ and the exponential connection in the dual coordinates $\theta^i = \xi(p_i) - \xi(p_n)$ \citep[Thm.~2.9 and eq.~(2.130)]{eguchikomori2022}. Thus every admissible $U$ remains dually flat. Since $C_U$ is affine in its first argument, its empirical average is $C_U(\bar x,q)$, the $U$-version of the two-stage reduction of \autoref{lem:two-stage} \citep[eqs.~(4.88)--(4.89)]{eguchikomori2022}.

The statistical structure enters through the family on which the
representation is affine.

\begin{definition}[$U$-model; {\citealp[Def.~3.1]{eguchikomori2022}}]\label{def:u-model}
For component scores $t: \{1,\dots,n\} \to \R^k$, the \emph{$U$-model}
generated by $t$ is
\begin{equation}\label{eq:u-model}
  E_t^U = \{p_\theta: \xi(p_\theta(i)) = \theta \cdot t(i) - b(\theta),
  \ i = 1,\dots,n\},
  \qquad
  p_\theta(i) = u\big(\theta \cdot t(i) - b(\theta)\big),
\end{equation}
the family whose $\xi$-representation is affine in the natural parameter
$\theta \in \R^k$, with $b(\theta)$ the unique value for which
$\sum_{i=1}^n u(\theta \cdot t(i) - b(\theta)) = 1$.
\end{definition}

The normalizer $b$ plays the role of the log-partition function, and at $U = \exp$ the family reduces to $p_\theta(i) \propto e^{\theta\cdot t(i)}$. Geometrically, $E_t^U$ is flat for the exponential-type connection because its $\xi$-representation is affine, while $M_\tau$ is flat for the mixture connection because it is affine in $p$ \citep[Thm.~3.3]{eguchikomori2022}. The generalized Pythagorean theorem governs the intersection of these two families.

\begin{theorem}[Generalized Pythagorean relation for $D_U$;
{\citealp[Thm.~2.10]{eguchikomori2022}}]\label{thm:pyth}
Let $p,q,r \in \simplex$. The $\nabla^{(m,U)}$-geodesic
from $p$ to $q$ is $g$-orthogonal at $q$ to the
$\nabla^{(e,U)}$-geodesic from $q$ to $r$ if and only if
\begin{equation}\label{eq:pyth-U}
  D_U[p:r] = D_U[p:q] + D_U[q:r].
\end{equation}
\end{theorem}

Write $\bar t = \E_{\bar x}[t]$ for the value of the statistic at $\bar x$. With $p = \bar x$, $r \in E_t^U$ and $q$ the member of $E_t^U$ with $\E_q[t] = \bar t$, the two geodesics are orthogonal at $q$ and \eqref{eq:pyth-U} gives $D_U[\bar x:r] \ge D_U[\bar x:q]$ with equality only at $r = q$ \citep[Thm.~3.4]{eguchikomori2022}: the projection is the point where the model meets the moment surface, the discrepancy to any other model member splitting into the part captured by the fit and the remaining within-model discrepancy. In the parameter, the same point solves the moment-matching equation
\begin{equation}\label{eq:u-momentenabgleich}
  \E_{p_{\hat\theta}}[t] = \bar t,
\end{equation}
since the divergence and the cross-entropy differ by $H_U(\bar x)$, free of $\theta$, and $\nabla_\theta C_U(\bar x,p_\theta) = \E_{p_\theta}[t] - \bar t$ \citep[Thm.~4.3]{eguchikomori2022}. Thus every admissible $U$ gives the same estimating equation; only the family $E_t^U$ in which it is solved depends on $U$. The maximizer of $H_U$ over $M_{\bar t}$ is the same point \citep[Thm.~3.1]{eguchikomori2022}: minimum $U$-divergence onto the model and maximum $U$-entropy over the moment surface give the same point $E_t^U \cap M_{\bar t}$, the estimator $\hat q_{\mathrm{EKO}} = p_{\hat\theta}$ fixed by \eqref{eq:u-momentenabgleich}.

This changes over the saturated (unrestricted) model, where the projection mechanism generates not one (as before) but two means. By \autoref{lem:two-stage}, $\arg\min_{q}\frac{1}{N}\sum_k D_U[x_k:q] = \bar x$ for every admissible $U$: the barycenter when minimizing the second argument is the arithmetic mean in mixture coordinates. Reversing the arguments, the minimizer $q^\ast$ instead averages the data in the dual coordinates and renormalizes, a quasi-arithmetic mean with generator $\xi$, unique in the interior by strict convexity; for Shannon, $\xi = \log$ and
\begin{equation}\label{eq:left-kl-barycenter}
  \nabla\psi_U(q^{\ast}) = \frac{1}{N}\sum_{k=1}^N\nabla\psi_U(x_k) + \mu\one
  \quad\Longrightarrow\quad
  q^{\ast}_i = \frac{\big(\prod_{k=1}^N (x_k)_i\big)^{1/N}}{\sum_{j=1}^n\big(\prod_{k=1}^N (x_k)_j\big)^{1/N}},
\end{equation}
the closed geometric mean of the compositions, the (bary)center of compositional data analysis \citep{pawlowskyglahnegozcue2002}, which \autoref{subsec:lcf} identifies with the Fr\'echet mean of the flat Aitchison-type geometry. The two barycenters thus average the data in the two flat coordinate systems of the dually flat structure, mixture and Legendre-dual respectively.

A well-known instance is the one-parameter family of generators of the power-law family \citep[Ex.~2, eqs.~(8)--(9)]{eguchikomoriohara2014}
\begin{equation}\label{eq:beta-family}
  U_\beta(s) = \frac{1}{1 + \beta}\,(1 + \beta s)_+^{(1+\beta)/\beta},
  \qquad
  u_\beta(s) = (1 + \beta s)_+^{1/\beta},
  \qquad
  \xi_\beta(t) = \frac{t^\beta - 1}{\beta},
  \qquad
  \beta \ge 0.
\end{equation}
The induced divergence $D_\beta$ is the density-power divergence of Basu, Harris, Hjort and Jones \citep[whose tuning
parameter is denoted $\alpha$; ours is $\beta$ to keep notation clear]{basu1998}, written out in \eqref{eq:dpd}, and the $U$-model is the power-law family $p_\theta(i) \propto (1 + \beta\,\theta \cdot t(i))_+^{1/\beta}$. Under \eqref{eq:u-divergenz} the generating potential is the separable $\psi_\beta(p) = \sum_i\big[p_i^{1+\beta}/(\beta(1 + \beta)) - p_i/\beta\big]$, with Hessian metric $g_{ij} = p_i^{\beta-1}\delta_{ij}$ interpolating between the Fisher-Rao metric at $\beta \to 0$ and the Euclidean metric at $\beta = 1$, where $D_1[p:q] = \tfrac{1}{2}\|p - q\|^2$. As $\beta \to 0$, $\xi_\beta \to \log$ and $u_\beta \to \exp$: $D_\beta \to D_{\mathrm{KL}}$, the model reverts to the exponential family, and
$\hat q_{\mathrm{EKO}}$ to the maximum-likelihood estimate.

For every $\beta$, however, the estimating equation over the \emph{matched} family, i.e.\ the $U$-model $E_t^{U_\beta}$ built from the same $U_\beta$ as the loss, so that loss and model share the representation $\xi_\beta$, still is the unweighted moment equation \eqref{eq:u-momentenabgleich}. This happens because powers cancel: the loss weights the model score by the $\beta$-th power of the cell probability (the weight that is the source of the robustness of the density-power divergence \citep{basu1998}), while on the matched family, where $p_\theta(i) = u_\beta(\theta \cdot t(i) - b(\theta))$ and $u_\beta' = u_\beta^{1-\beta}$, the score itself carries the reciprocal power:
\begin{equation}\label{eq:beta-kuerzung}
  p_\theta(i)^{\beta}\,\nabla_\theta\log p_\theta(i)
  = p_\theta(i)^{\beta}\,p_\theta(i)^{-\beta}\,
   \big(t(i) - \nabla b(\theta)\big)
  = t(i) - \nabla b(\theta).
\end{equation}
Robustness is therefore a property of the loss-family pair; \autoref{subsec:robustness} shows how it reappears when the $\beta$-loss is decoupled from its matched family and fitted instead to the undeformed exponential family $E_t$, where no such compensation occurs.

In the Shannon case $\beta = 0$, $p_\theta(i) = e^{\theta\cdot t(i)}/Z(\theta)$ and stationarity of the log-likelihood $\theta \cdot \bar t - \log Z(\theta)$ is $\E_{p_{\hat\theta}}[t] = \bar t$. Throughout we illustrate with $n = 3$ and $t(i) = i$, the category label itself, and the interior sample $x_1 = (0.7,0.1,0.2)$, $x_2 = (0.4,0.3,0.3)$, $x_3 = (0.4,0.2,0.4)$, of mean $\bar x = (0.5,0.2,0.3)$ and hence $\bar t = 1.8$, a three-category example; this gives the maximum-entropy distribution of that mean, the first row of \autoref{tab:fits} (\eqref{eq:eko-numerisch} in \autoref{app:eko-shannon-example}).

The same identity has a second consequence: it closes the dually flat case under Bayesian updating. Loss-based updating \citep{bissiriholmeswalker2016} multiplies a belief distribution over $q$ by the exponentiated negative empirical risk, and reduces to Bayes' theorem exactly when the loss is the logarithmic one of \autoref{subsubsec:likelihood}. For an admissible $U$, an anchor $m \in \simplex$ and $\nu > 0$ in the role of a prior sample size, set
\begin{equation}\label{eq:u-gibbs}
  \operatorname{prior}_{\nu,m}(q) \propto \exp\big\{-\nu\,D_U[m:q]\big\},
  \qquad
  \operatorname{post}(q) \propto \operatorname{prior}_{\nu,m}(q)\,
  \exp\Big\{-\sum_{k=1}^N D_U[x_k:q]\Big\},
\end{equation}
whenever the normalization is finite.

\begin{proposition}[Conjugacy of the forward $U$-updating]\label{prop:u-conjugacy}
The family \eqref{eq:u-gibbs} is closed under its own updating,
\begin{equation}\label{eq:u-conjugacy}
  \operatorname{post} = \operatorname{prior}_{\nu',m'},
  \qquad \nu' = \nu + N,
  \qquad m' = \frac{\nu\,m + N\,\bar x}{\nu + N},
\end{equation}
the posterior anchor being the mixture-geodesic combination of the prior anchor and the sample mean, the prior entering as $\nu$ pseudo-observations placed at $m$. Each $\operatorname{prior}_{\nu,m}$ has mode $m$, so the maximum-a-posteriori estimate over a model $M$ is the divergence projection $\arg\min_{q\in M}D_U[m':q]$, the estimator of this subsection applied to the prior-blended mean. At $U = \exp$ and vertex data the updating is exact Bayes, $\operatorname{prior}_{\nu,m}$ being the Dirichlet density at $\alpha = \nu m + \one$.
\end{proposition}

Conjugacy is therefore not a Dirichlet accident but a property of the principle itself: apply \autoref{lem:two-stage} to the points $(m,x_1,\dots,x_N)$ with weights $(\nu,1,\dots,1)$. It stops at the logarithmic divergences, which do not factor through the mean \citep{banerjeeguowang2005}.

This estimation principle is built on the dual flatness of $D_U$ on the model. The next one keeps the projection but deforms the duality itself, and can return a different estimate outside the undeformed limit.

\subsection{Projectively flat estimation: the WZTY principle}\label{subsec:wzty}

The deformation functions
\begin{equation}\label{eq:log-exp-lambda}
  \log_\lambda(u) = \frac{u^{\lambda} - 1}{\lambda},
  \qquad
  \exp_\lambda(s) = (1 + \lambda s)_+^{1/\lambda},
  \qquad
  \lambda > 0,
\end{equation}
mutually inverse and reducing to $(\log,\exp)$ as $\lambda \to 0$, are the pair $(\xi_\beta,u_\beta)$ of \eqref{eq:beta-family} at $\beta = \lambda$.
The estimation principle of this subsection is nevertheless not an instance of the previous one; the difference lies in where the functions \eqref{eq:log-exp-lambda} enter. In the $U$-construction $\xi_\lambda$ deforms the dual coordinate: the Legendre dual of the mixture coordinate $p_i$ is $\xi_\lambda(p_i)$, the model is affine in it, and the pairing between primal and dual variables remains the bilinear $x \cdot y$ that underlies every Bregman divergence. Here, following Wong and Zhang \citep{zhangwong2022}, the probabilities enter the divergence untransformed and the deformation happens in the pairing, which becomes the logarithmic cost $c^{(\lambda)}(x,y) = \tfrac{1}{\lambda}\log(1 + \lambda\,x \cdot y) = \log\exp_\lambda(x \cdot y)$, the divergence being the logarithmic one of \autoref{subsubsec:metric-lconvex}; in the robust family of \autoref{subsec:eko} the exponential tilt $e^{s}$ becomes the power-law tilt $\exp_\lambda(s)$. The geometry beneath the projection is projectively rather than dually flat, in the terminology of Pal and Wong.

\begin{definition}[$\lambda$-exponential family]\label{def:lambda-expfam}
Replacing $\exp$ by $\exp_\lambda$ in the exponential family of \autoref{subsec:eko}, the \emph{$\lambda$-exponential family} generated by the component scores $t: \{1,\dots,n\} \to \R^k$ is
\begin{equation}\label{eq:lambda-expfamilie}
  E_t^{\lambda} = \{p_\theta: p_\theta(i) =
  \exp_\lambda\big(\theta \cdot t(i) - \varphi_\lambda(\theta)\big)\},
\end{equation}
with $\varphi_\lambda(\theta)$ the unique value for which $\sum_{i=1}^n\exp_\lambda\big(\theta \cdot t(i) - \varphi_\lambda(\theta)\big) = 1$:
the (subtractive)\footnote{We use the subtractive normalization throughout; the divisive variant normalizes by a factor instead, $p_\theta(i) = \exp_\lambda(\theta \cdot t(i))/Z_\lambda(\theta)$; the relation between the two forms is explored in \citep{zhangwong2022}.\label{fn:subtractive}} normalizer of the family, the analogue of $b(\theta)$ in the $U$-model \eqref{eq:u-model}, reducing to the log-partition function $\log Z(\theta)$ as $\lambda \to 0$, i.e. as the family tends to the exponential one.
\end{definition}

The canonical divergence of the family is the $\lambda$-logarithmic divergence \citep{palwong2018,wong2018},
\begin{equation}\label{eq:lambda-ldiv}
  D^{(\lambda)}_\phi(p \mid q) = \frac{1}{\lambda}\log\big(1 + \lambda\,
  \nabla\phi(q) \cdot (p - q)\big) - \big(\phi(p) - \phi(q)\big),
\end{equation}
defined for a $\lambda$-exponentially concave generator $\phi$, i.e. one with $e^{\lambda\phi}$ concave: the $L$-divergence \eqref{eq:l-divergence} at $\lambda = 1$ and, since $\tfrac{1}{\lambda}\log(1 + \lambda s) \to s$, the Bregman divergence \eqref{eq:bregman} of $-\phi$ as $\lambda \to 0$.

The Legendre duality of \autoref{subsec:eko} is conjugacy against the bilinear cost $c_0(\theta,\eta) = \theta \cdot \eta$, natural parameter paired with dual variable, with the Bregman divergence as the Fenchel-Young remainder. In the present setting it becomes the $c$-conjugacy of optimal transport against the logarithmic cost $c^{(\lambda)}$, on $1 + \lambda\,x \cdot y > 0$; the $c$-transform, the $\lambda$-Legendre map and the representation of $D^{(\lambda)}_\phi$ as the corresponding remainder are developed by Wong \citep{wong2018}, each term tending to its Legendre counterpart as $\lambda \to 0$.

The first consequence of the duality is the analogue of the Pythagorean relation of \autoref{thm:pyth}.

\begin{theorem}[$\lambda$-Pythagorean relation;
\citealp{palwong2018,wong2018}]\label{thm:lambda-pyth}
Let $p,q,r \in \simplex$. If the primal geodesic from $p$ to $q$ is $g$-orthogonal at $q$ to the dual geodesic from $q$ to $r$, then
\begin{equation}\label{eq:lambda-pyth}
  D^{(\lambda)}_\phi(p \mid r) = D^{(\lambda)}_\phi(p \mid q)
  + D^{(\lambda)}_\phi(q \mid r).
\end{equation}
\end{theorem}

For $\lambda = 1$ this is the Pythagorean theorem of Pal and Wong; the general case is due to Wong. As in \autoref{thm:pyth}, the relation makes the projection the unique minimizer along the relevant geodesics. The geometry beneath that mechanism changes: the computation of \autoref{app:ldiv-eguchi} with $A = 1 + \lambda\,\nabla\phi(q) \cdot (p - q)$ and $\Phi = e^{\lambda\phi}$ gives
\begin{equation}\label{eq:lambda-metric}
  g = g^{\mathrm H} - \lambda\,\nabla\phi\,\nabla\phi\T
  = -\tfrac{1}{\lambda}\,e^{-\lambda\phi}\,\nabla^{2}e^{\lambda\phi},
\end{equation}
the rank-one correction of \eqref{eq:ldiv-metric} with weight $\lambda$: positive semidefiniteness is $\lambda$-exponential concavity restated, and the correction vanishes into the Bregman Hessian as $\lambda \to 0$. The primal connection is no longer flat; instead:

\begin{proposition}[Projective flatness of the $\lambda$-logarithmic connection; \citealp{wong2018,palwong2018}]\label{prop:proj-flat}
For every $\lambda$-exponentially concave generator $\phi$, the mixture connection of the $\lambda$-logarithmic divergence \eqref{eq:lambda-ldiv}
satisfies, in the coordinates $p$ in which the Bregman mixture connection vanishes,
\begin{equation}\label{eq:christoffel-proj}
  \Gamma^{(m)\,l}_{ij} = \delta^l_i\,a_j + \delta^l_j\,a_i,
  \qquad
  a = -\lambda\,\nabla\phi,
\end{equation}
and its curvature tensor is
\begin{equation}\label{eq:ldiv-curvature}
  R^{(m)\,l}{}_{kij} = \lambda\,\big(\delta^l_j\,g_{ik} - \delta^l_i\,g_{jk}\big),
\end{equation}
the constant-curvature form $R^{l}{}_{kij} = \kappa\,(\delta^l_i g_{jk} - \delta^l_j g_{ik})$ with $\kappa = -\lambda$, the
Pal-Wong value $\kappa = -1$ at $\lambda = 1$.\footnote{Contracting the difference tensor $\delta^l_ia_j + \delta^l_ja_i$ with $\dot\gamma^i\dot\gamma^j$ yields $2\,(a \cdot \dot\gamma)\,\dot\gamma^l$, parallel to the velocity. This is Weyl's characterization of when two symmetric connections share their geodesics.}.
\end{proposition}

Both displays are obtained directly from \eqref{eq:lambda-ldiv}, by differentiating at the diagonal and substituting the metric relation $\phi_{ij} = -g_{ij} - \lambda\,\phi_i\phi_j$ of \eqref{eq:lambda-metric}, which cancels the cubic terms, for every $n$ and every $\lambda > 0$ (\autoref{app:ldiv-projective}). The modification \eqref{eq:christoffel-proj} means that the connection shares its geodesics (see the footnote) with the flat Bregman connection as unparametrized curves: straight lines in the $p$-chart, traversed at the non-affine rate $dt \propto e^{-2\lambda\phi}d\sigma$ in the Euclidean parameter $\sigma$ (\autoref{app:wzty-reparam}). Thus the structure is dually projectively flat with constant dual curvature $-\lambda$; as $\lambda \to 0$ both the projective term and the curvature vanish, recovering the dually flat case. EKO and WZTY therefore share the additive Pythagorean projection but differ in the geometry: linear versus $c^{(\lambda)}$-duality, bare-Hessian versus conformal-Hessian metric, affine versus projective flatness.

The estimator is the projection of the sample onto the deformed family,
\begin{equation}\label{eq:wzty-schaetzer}
  \hat q_{\mathrm{WZTY}} = p_{\hat\theta},
  \qquad
  \hat\theta = \arg\min_{\theta}\,\frac{1}{N}\sum_{k=1}^N
  D^{(\lambda)}_\phi(x_k \mid p_\theta),
\end{equation}
whose uniqueness along the relevant geodesics follows from \autoref{thm:lambda-pyth}; the full projection theory onto primal and dual autoparallel submanifolds is developed by Tao and Wong \citep{taowong2021}. Unlike in the dually flat case, the generator is part of the estimator: the target depends on the pair $(\phi,E_t^{\lambda})$, a point taken up in \autoref{subsec:robustness}.

Differentiating \eqref{eq:lambda-ldiv} in its second argument, $\phi(q)$ and $\log A_k$ depend on $q$, where $A_k = 1 + \lambda\,\nabla\phi(q) \cdot (x_k - q)$; the derivative of the logarithm puts the reciprocal $w_k(q) := A_k^{-1}$ in front of each observation, and the identity $1 - w_k = \lambda\,w_k\,\nabla\phi(q) \cdot (x_k - q)$ turns the gradient into $\nabla_qD^{(\lambda)}_\phi(x_k \mid q) = -w_k(q)\,g(q)\,(x_k - q)$ with $g$ the metric \eqref{eq:lambda-metric} (\autoref{app:wzty-score}). Stationarity in the directions tangent to the family, with $J_\theta = \partial_\theta p_\theta$ the $n \times k$ Jacobian of the parametrization, whose columns are tangent to the simplex, is therefore
\begin{equation}\label{eq:lambda-score}
  J_\theta\T\sum_{k=1}^N w_k(q)\,g(q)\,(x_k - q) = 0
  \quad\text{at }q = p_{\hat\theta} = \hat q_{\mathrm{WZTY}},
\end{equation}
and for the saturated model this is stationarity in every simplex tangent direction: weights and metric are the same deformation, and the score sets the reweighted residuals $g$-orthogonal to the family. The sample enters the score only through the self-normalized reweighted mean $\sum_k\hat w_k(q)\,x_k$, $\hat w_k = w_k/\sum_lw_l$, a reweighting the estimate itself selects (\autoref{app:wzty-score}). As $\lambda \to 0$ every $w_k \to 1$, the reweighted mean collapses to $\bar x$, the divergence tends to the Bregman divergence of $-\phi$, and the projection of the sample reduces, by \autoref{lem:two-stage}, to the projection of $\bar x$: the estimator is again the dually flat one.

At $\lambda = 1$ the subtractively-normalized (footnote~\ref{fn:subtractive}) family is affine,
\begin{equation}\label{eq:lambda1-family}
  p_\theta(i) = \tfrac{1}{n} + \theta\,(t(i) - \bar t_0),
  \qquad \bar t_0 = \tfrac{1}{n}\sum_{j=1}^n t(j).
\end{equation}
For $n = 3$, $t(i) = i$ and the uniform geometric-mean generator, the running sample gives the WZTY row of \autoref{tab:fits}, which differs from $\hat q_{\mathrm{EKO}}$; the minimization and induced weights are in \autoref{app:wzty-interior-example}.
For categorical data this projection estimator is not well-defined for the generators we introduced: $\phi_\pi(e_i) = -\infty$, hence $D^{(\lambda)}_{\phi_\pi}(e_i \mid q) = +\infty$. One could drop the data-side infinite term to produce a finite penalized objective but that would not give the likelihood of the counts; the categorical version developed by Wong, Zhang, Tian and Yang is instead a maximum likelihood over the deformed family,
\begin{equation}\label{eq:lambda-mle}
  \hat q_{\lambda\mathrm{ML}} = p_{\hat\theta},
  \qquad
  \hat\theta = \arg\max_{\theta}\,\sum_{i=1}^n\hat p_i\log p_\theta(i),
\end{equation}
Its score equation is not the raw moment matching \eqref{eq:u-momentenabgleich} that maximum likelihood produces over the undeformed family: the pair is decoupled in the sense of \autoref{subsec:eko}, (the $\lambda = 0$ loss fitted to the $\lambda$-deformed family), so the factor $p_{\hat\theta}^{-\lambda}$ coming from the representation does not cancel as on the matched pair and is still present in the score (\autoref{app:wzty-escort-derivative}), and the first-order condition of the $\lambda$-likelihood maximization \eqref{eq:lambda-mle} is the escort-moment equation
\begin{equation}\label{eq:begleit-moment}
  \frac{\sum_i \hat p_i\,p_{\hat\theta}(i)^{-\lambda}\,t(i)}{\sum_i \hat p_i\,p_{\hat\theta}(i)^{-\lambda}}
  = \nabla\varphi_\lambda(\hat\theta),
\end{equation}
matching a $p_{\hat\theta}^{-\lambda}$-escort mean of $t$ to the model's $\lambda$-mean rather than to $\bar t$. As $\lambda \to 0$ the weights tend to one and \eqref{eq:begleit-moment} returns to \eqref{eq:u-momentenabgleich}; the same escort structure, with fixed-point algorithms and continuous analogues, appears in \citep{tian2026,guilmeau2024}\footnote{The term ``escort'' comes from statistical mechanics \citep{beckschlogl1993}: the tilted distributions $p^{(a)}_i = p_i^{a}/\sum_jp_j^{a}$ accompany (escort) $p$, shifting weight to high-probability cells for $a > 1$ and to low-probability cells as $a$ decreases; in \eqref{eq:begleit-moment} the order is $-\lambda$, so the escort mean weights the cells the fitted model sees as unlikely.}.

For the running example $\lambda = 1$, the categorical likelihood version solves the escort equation at the $\lambda$-ML row of \autoref{tab:fits}, with $\E_{\hat q_{\lambda\mathrm{ML}}}[t] = 11/6$; the short scalar calculation is in \autoref{app:wzty-categorical-example}. Thus likelihood and divergence projection give different estimates even when the categorical and interior samples agree in mean: the two $\lambda$-estimators are not equivalent.

More generally, for weights $\pi \in \simplex$ and order $\beta \le 1$ the power mean $M^\pi_\beta(x) = (\sum_i\pi_ix_i^\beta)^{1/\beta}$ is concave, with $\beta = 0$ understood by continuity, so $\phi = \log M^\pi_\beta$ is an exponentially concave generator.

The deformation widens the choice of generator in a second, independent direction, and the two directions compose: the order $\beta$ and the deformation $\lambda$ have so far moved together ($(\xi_\beta,u_\beta)$ at $\beta = \lambda$ is $(\log_\lambda,\exp_\lambda)$) or not at all; letting them move independently, the density-power potential of \autoref{subsec:eko} becomes itself a candidate generator for the $\lambda$-divergence.

\begin{proposition}[The density-power generators parametrized by $\lambda$]\label{prop:beta-lambda}
Fix $\beta > 0$, $\lambda > 0$, and let $\phi_\beta := -\psi_\beta$, the sign-flipped density-power potential of \autoref{subsec:eko}, with $\xi_\beta$ as in \eqref{eq:beta-family}. Then $\phi_\beta$ is $\lambda$-exponentially concave on $\simplex$ if and only if
\begin{equation}\label{eq:beta-grenze}
  \lambda\,\sup_{q\in\simplex}W_\beta(q) \le 1,
  \qquad
  W_\beta(q) := \min_{c\in\R}\sum_{i=1}^{n}q_i^{1-\beta}
  \big(\xi_\beta(q_i) - c\big)^{2},
\end{equation}
the constant being sharp, since $W_\beta(q)$ is attained by a tangent direction at every $q$; and the induced metric \eqref{eq:lambda-metric} is
\begin{equation}\label{eq:beta-lambda-metrik}
  g^{(\beta,\lambda)}_{ij}(q)
  = q_i^{\beta-1}\,\delta_{ij} - \lambda\,\xi_\beta(q_i)\,\xi_\beta(q_j),
\end{equation}
the $U_\beta$-metric \eqref{eq:u-metrik} as bare-Hessian part, with $\lambda$ entering only through the rank-one correction.
\end{proposition}

By \autoref{prop:proj-flat} the structure is dually projectively flat with dual curvature $-\lambda$ for every $\beta$, and the divergence has a closed form in the escort weights $q_i^{1+\beta}/\sum_jq_j^{1+\beta}$ and the ratios $p_i/q_i$ (\eqref{eq:beta-lambda-div} in \autoref{app:est-proofs}). Its two limits are the two principles of this section: as $\lambda \to 0$ the density-power divergence $D_\beta$, as $\beta \to 0$ a $\lambda$-tempered Kullback-Leibler divergence \eqref{eq:temperiertes-kl}, itself $D_{\mathrm{KL}}$ as $\lambda \to 0$. At that edge \eqref{eq:beta-grenze} reads $\lambda\,\sup_q\Var_q(\log q) \le 1$, which at $\lambda = 1$ holds exactly for $n \le 3$; at $\beta = 1$ it reads $\lambda \le \tfrac{n}{n - 1}$, so the Euclidean generator is admissible at the Pal-Wong point for every $n$ (values for $n = 4,5$ in \autoref{app:est-proofs}).

The weights $w_k$ have a second reading, on the prior side. The updating rule \eqref{eq:u-gibbs} applies unchanged in this case, but its closure fails. With a flat prior,
\begin{equation}\label{eq:lambda-gibbs}
  \operatorname{post}(q)\;\propto\;
  e^{-N\phi(q)}\prod_{k=1}^N A_k(q)^{-1/\lambda},
  \qquad
  A_k(q) = 1 + \lambda\,\nabla\phi(q) \cdot (x_k - q),
\end{equation}
each data factor the reciprocal $\lambda$-exponential of a linear residual: polynomial, Student-type tails, against the exponential Bregman factors. Its mode is \eqref{eq:wzty-schaetzer} and its data factors are the score weights $w_k = A_k^{-1}$ themselves; for the running generators the posterior is proper whenever $\lambda\pi_i < 1$ (\autoref{app:lambda-bayes-propriety}).

Closure returns in a randomized sense. Since $A$ is affine in its data argument, $\sum_ku_kA_k(q)$ is the value $A(\bar x_u,q)$ at the reweighted mean $\bar x_u = \sum_ku_kx_k$, and the Feynman parametrization of a product of powers \citep{feynman1949} turns \eqref{eq:lambda-gibbs} into
\begin{equation}\label{eq:randomisierte-konjugiertheit}
  \operatorname{post}(q)\;\propto\;
  \E_{u\sim\Dir(1/\lambda,\dots,1/\lambda)}
  \Big[\,e^{-N\phi(q)}\,A(\bar x_u,q)^{-N/\lambda}\,\Big],
\end{equation}
an exact mixture of conjugate kernels anchored at Dirichlet-reweighted sample means: draw weights, update conjugately, average. As $\lambda \to 0$ the weight law concentrates at the uniform weights and \eqref{eq:randomisierte-konjugiertheit} returns to \autoref{prop:u-conjugacy}; at $\lambda = 1$ it is the uniform law on the weight simplex, the Bayesian bootstrap of Rubin \citep{rubin1981} (numerical check in \autoref{app:lambda-bayes-numerics}).

\subsection{The Levi-Civita-Fr\'echet principle}\label{subsec:lcf}

EKO and WZTY use divergence losses; the Levi-Civita-Fr\'echet principle takes instead the squared geodesic distance $\ell(q,x) = d(q,x)^2$, giving the estimator
\begin{equation}\label{eq:lcf-estimator}
\hat q_{\mathrm{LCF}} = \arg\min_{q\in M}\frac{1}{N}\sum_{k=1}^N d(q,x_k)^2,
\end{equation}
the Fr\'echet mean of the sample in the metric $g$. With a real covariate the same loss defines geodesic regression \citep{fletcher2013}, the least-squares fit of a Levi-Civita geodesic $s \mapsto \exp_q(s\,v)$ to responses on a Riemannian manifold, whose intercept-only case $v = 0$ is \eqref{eq:lcf-estimator}. Fletcher fixes the metric and lets the manifold be general; here the manifold is the simplex and the metric is what varies, so the covariate is left out and the comparison runs across geometries; regression on the simplex with a covariate is treated in \citep{navasmest}. Differentiation of $d^2$ is Riemannian: the exponential map replaces the global Legendre coordinates.\footnote{In the self-dual structure $(g,\nabla^{(\mathrm{LC})})$, $\tfrac{1}{2}d^2$ is the canonical divergence of Ay and Amari \citep{ayamari2015}; the corresponding characterization theory is developed in \citep{navasmest}.}

For the differentiation we write $\log_q: \mathcal U_q \to T_q\simplex$, $\log_q(x) = \dot\gamma(0)$, for the minimizing geodesic $\gamma: [0,1] \to \simplex$ of $\nabla^{(\mathrm{LC})}$ with $\gamma(0) = q$, $\gamma(1) = x$, defined on the open set $\mathcal U_q$ of points reached from $q$ by a unique such geodesic.

\begin{proposition}[Gradient of squared distance]\label{prop:grad-dist}
For $x \in \simplex$ and any $q$ with $x \in \mathcal U_q$,
\begin{equation}\label{eq:grad-dist-sq}
  \nabla_q\, d(q,x)^2 = -2\log_q(x).
\end{equation}
\end{proposition}

Summing over the data gives the stationarity condition of the unconstrained problem, $\frac{1}{N}\sum_{k=1}^N \log_{\hat q_{\mathrm{LCF}}}(x_k) = 0$, the geodesic-loss analogue of EKO's moment matching and of WZTY's reweighted-mean condition \eqref{eq:lambda-score}: each is the score equation of the same M-estimator \eqref{eq:m-estimator}, differing only in the loss differentiated. For a proper submanifold $M$ the condition holds in the directions tangent to $M$ at the estimate. The cases below are worked directly, without computing $\log_q$ itself.

Existence depends on how the chosen geometry treats the boundary, and one criterion covers all generators. For a boundary point $e$ of the simplex, let $d(q,e)$ be the infimum of lengths of interior paths from $q$ converging to $e$; each $d(\cdot,e)$ is $1$-Lipschitz, extends continuously to the metric completion, and its finiteness does not depend on $q$, by the triangle inequality.

\begin{proposition}[Existence dichotomy]\label{prop:existence-dichotomy}
Let the data $y_1,\dots,y_N$ lie in the closed simplex, distances to boundary points understood as above. If some $y_k$ lies at infinite distance, the objective is identically $+\infty$ and no estimator exists. If every $y_k$ lies at finite distance and the metric completion of $(\simplex,g)$ is compact, the objective extends continuously to the completion and attains its minimum there.
\end{proposition}

For interior data the first alternative is vacuous (\autoref{app:rob-existence-proof}); which alternative occurs depends on the generator.

(a) For Fisher-Rao, the isometry $\Pi(p) = 2\sqrt{p}$ onto the sphere of radius $2$ gives $\langle\Pi(p),\Pi(q)\rangle = 4\sum_i\sqrt{p_iq_i}$, hence
\begin{equation}\label{eq:rao-abstand}
  d(p,q) = 2\arccos\Big(\sum_{i=1}^n\sqrt{p_iq_i}\Big),
  \qquad
  d(p,e_i) = 2\arccos(\sqrt{p_i}) < \pi,
\end{equation}
the bound \eqref{eq:diameter} being attained, since $d(e_i,e_j) = \pi$ for $i \neq j$, while every interior point stays strictly inside it. On the saturated model, differentiating the Fisher-Rao objective under the simplex constraint gives (\autoref{app:lcf-fr-stationarity})
\begin{equation}\label{eq:lcf-stationarity-fr}
  \frac{1}{N}\sum_{k=1}^N\frac{\theta_k(q)}{\sin\theta_k(q)}
  \sqrt{\frac{(x_k)_i}{q_i}} = \nu,
  \quad i = 1,\dots,n,
  \qquad
  \theta_k(q) = \arccos\Big(\sum_{j=1}^n\sqrt{q_j(x_k)_j}\Big),
\end{equation}
transcendental for $n \ge 3$: there is no closed form, and the estimator is found by solving a system of equations\footnote{The condition is $\operatorname{grad}_qF = 0$ in $T_q\simplex$ for $F(q) = \frac{1}{N}\sum_kd(q,x_k)^2$; \eqref{eq:lcf-stationarity-fr} writes it in the ambient coordinates of $\R^n$, with $\nu$ the multiplier of $\sum_iq_i = 1$, so the $n$ equations holding with one common $\nu$ is the vanishing of the tangential part. The restriction $x \in \mathcal U_q$ of \autoref{prop:grad-dist} is vacuous, the injectivity radius of the ambient sphere, $2\pi$, exceeding the diameter $\pi$.}. Vertices lie at finite distance by \eqref{eq:rao-abstand}, so the completion is the closed orthant of the sphere, compact, and the second alternative of \autoref{prop:existence-dichotomy} is the one occurring; the minimizer is interior for interior data (\autoref{app:lcf-fr-stationarity}). Uniqueness needs no smallness condition on the sample: every point of the positive orthant lies within $2\arccos(1/\sqrt{n}) < \pi$ of the orthant's center, inside the uniqueness ball of Afsari \citep[Thm.~2.1]{afsari2011} (earlier, smaller balls in \citealp{karcher1977,kendall1990}).

By standard convexity arguments, on a space of constant curvature $\kappa > 0$, the Hessian of $d(\cdot,x)^2$ at distance $d$ from $x$ has eigenvalue
\begin{equation}\label{eq:sphere-hessian-ev}
  2\ \text{ along the geodesic to }x,
  \qquad
  2\sqrt{\kappa}\,d\cot(\sqrt{\kappa}\,d)\ \text{ on its orthogonal complement},
\end{equation}
both positive while $\sqrt{\kappa}\,d < \pi/2$ (\autoref{app:sphere-hessian}); a positively weighted sum of such functions is strictly geodesically convex and has at most one stationary point. Thus every interior sample has a unique LCF mean; full-support categorical samples are covered by the vertex distance in \eqref{eq:rao-abstand}, and zero-frequency cells reduce the problem to the corresponding face.

On the running example the Fisher-Rao LCF mean differs from $\bar x$, and the categorical sample with the same mean gives a different mean again (\autoref{tab:fits}; \eqref{eq:lcf-numeric-fr}--\eqref{eq:lcf-numeric-cat} in \autoref{app:lcf-running-examples}): curvature makes the saturated Fr\'echet mean depend on more than the arithmetic mean.

(b) The second case is the weighted-geometric-mean generator, whose metric \eqref{eq:ldiv-beispiel-metrik} is flat by the pullback isometry $g = \clr_\pi^{\ast}\langle\cdot,\cdot\rangle_\pi$ onto $V_\pi$ \eqref{eq:ldiv-clr-rueckzug}.

\begin{proposition}[Fr\'echet mean under a flat pullback metric]\label{prop:flat-frechet}
Let $g = F^{\ast}\langle\cdot,\cdot\rangle$ for a diffeomorphism $F$ from $\simplex$ onto an inner-product space $(V,\langle\cdot,\cdot\rangle)$. For interior points $x_1,\dots,x_m$ and weights $w_1,\dots,w_m \ge 0$ summing to one, the unconstrained Fr\'echet mean exists, is unique, and is the $F$-preimage of the $F$-average of the data, \eqref{eq:lcf-coda}.
\end{proposition}

Flatness gives the closed form that curvature prevented in (a) (\autoref{app:lcf-flat-pullback}):
\begin{equation}\label{eq:lcf-coda}
  q^{\ast} = F^{-1}\Big(\sum_{i=1}^m w_i\, F(x_i)\Big),
  \quad
  \hat q_{\mathrm{LCF}} = \clr_\pi^{-1}\Big(\tfrac{1}{N}\sum_{k=1}^N\clr_\pi(x_k)\Big),
\end{equation}
the second with $F = \clr_\pi$, which for uniform weights is the closed geometric mean \eqref{eq:left-kl-barycenter}, the CoDA center \citep{pawlowskyglahnegozcue2002} and the left Kullback-Leibler barycenter, these three constructions coinciding in the flat geometry (\autoref{tab:fits}).

The hypothesis fails at the boundary: every vertex lies at infinite $g$-distance from the interior as $p \to e_k$, $\clr_\pi(p)_k = (1 - \pi_k)\log p_k - \sum_{j\ne k}\pi_j\log p_j \to +\infty$ and a sample with a zero part falls under the first alternative of \autoref{prop:existence-dichotomy}. This is the geometric form of the CoDA zero problem \citep{martinfernandez2003}, and the contrast with Fisher-Rao is the dichotomy itself: the same data are fitted in one geometry and have no estimator in the other.

(c) Between the two previous cases lies the square-root generator $\phi(p) = \log(\sum_i\sqrt{p_i})$, half the logarithm of the uniform order-$\tfrac{1}{2}$ power mean up to an additive constant, which reaches its boundary at finite distance as Fisher-Rao does, but whose curvature is not constant and changes sign for $n = 3$. Its metric and fourth-root chart are
\begin{equation}\label{eq:wurzel-metrik}
g_{ij}(p) = \frac{\delta_{ij}}{4\rho p_i^{3/2}},
\qquad \rho = \sum_i\sqrt{p_i},
\qquad
s(p) = (p_i^{1/4})_i \in S = \{s \in \R^n: \textstyle\sum_i s_i^4 = 1\}.
\end{equation}

\begin{proposition}[Compactification of the square-root geometry]
\label{prop:sqrt-sphere}
The fourth-root chart $s$ is an isometry of $(\simplex,g)$ onto the open positive piece $S_+ = S \cap \R^n_{>0}$ of the compact hypersurface $S$ with the metric $\bar g_s = 4\|s\|^{-2}\langle\cdot,\cdot\rangle|_{T_sS}$, smooth on all of $S$; every minimizing geodesic of $(S,\bar g)$ between points of $S_+$ stays in $S_+$, so the intrinsic distance of $(\simplex,g)$ is that of $S$ and is realized by interior geodesics; and the metric completion of $(\simplex,g)$ is the compact piece $\bar S_+ = S \cap \R^n_{\ge0}$, to which the curvature extends smoothly, with positive injectivity radius.
\end{proposition}

The second alternative of \autoref{prop:existence-dichotomy} therefore holds here too, now with interior minimizing geodesics and a positive injectivity radius, which is what the local theory of \autoref{subsec:robustness} needs of a curved generator (proof and the constants for $n = 3$ in \autoref{app:est-proofs} and \autoref{app:rob-sqrt-certificates}).

Computationally, Fisher-Rao has a closed distance and a globally geodesically convex objective, the flat case has the closed form of \autoref{prop:flat-frechet}, and a general generator needs one two-point geodesic boundary-value problem per observation and objective call, whose solvability \autoref{prop:sqrt-sphere} guarantees for the square-root geometry and which requires additional assumptions in general.

\subsection{The symmetrized cost as a loss: the $D^{\mathrm{sym}}$-barycenter}\label{subsec:jbary}

The symmetrized cost \eqref{eq:symkosten} supports neither projection principle (\autoref{prop:symcost-exclusion}), while its induced metric is a rescaling of the flat case (b) of \autoref{subsec:lcf}; an LCF estimator would therefore return the same geometric mean \eqref{eq:lcf-coda}; what the cost adds to its quadratic part, its nonquadratic form with Hilbert-slope tails, enters only when the loss itself is minimized:
\begin{equation}\label{eq:jbary}
  \hat q_{\mathrm{sym}}
  = \arg\min_{q\in M}\frac{1}{N}\sum_{k=1}^N D^{\mathrm{sym}}(q,x_k),
\end{equation}
the $D^{\mathrm{sym}}$-barycenter, a Fr\'echet-type principle applied to the divergence rather than to the squared distance it agrees with to second order.

\begin{proposition}[Existence and uniqueness of the barycenter]\label{prop:jbary}
For interior observations $x_1,\dots,x_N$, the estimator \eqref{eq:jbary} over the saturated model exists, is unique, and is interior. The same holds over every model that is affine in the clr chart, provided the parametrization is injective; in particular over every exponential family $E_t$, whose clr representation $\clr(p_\theta)(i) = \theta \cdot \big(t(i) - \bar t_0\big)$, $\bar t_0 = \frac{1}{n}\sum_jt(j)$, is linear in $\theta$, the log-partition function cancelling under centering.
\end{proposition}

No Pythagorean identity is needed, as on log-affine families the convexity is global. The score is, by \eqref{eq:symkosten-clr},
\begin{equation}\label{eq:jbary-score}
  \nabla_\xi F(\xi)
  = \frac{1}{N}\sum_{k=1}^N
  \Big(\clr^{-1}(\xi - \xi_k) - \clr^{-1}(\xi_k - \xi)\Big),
\end{equation}
a difference of two compositions for each observation, hence of Euclidean norm at most $\sqrt{2}$, the value at two distinct vertices, so the loss is globally Lipschitz in the chart where the case (b) objective is quadratic. Together with the quadratic part \eqref{eq:symcost-germ}, this is a Huber-type loss whose crossover from quadratic near the center to linear in the tails is fixed by the cost. Hence concentrated data recover the CoDA/left-KL center to fourth order, while remote observations act through direction rather than magnitude. On the running sample the barycenter lies at clr-distance $0.006$ from the geometric mean (\autoref{tab:fits}; \eqref{eq:jbary-numeric} in \autoref{app:jbary-running-example}); under contamination toward a vertex the two estimators separate, as quantified in \autoref{subsubsec:rob-sym}. At the population level, the barycenter is a location functional of its own and, because the score is an odd function, it agrees with the CoDA center for clr laws symmetric about their center.

The barycenter is equivariant under perturbation, permutation and inversion, by \autoref{subsubsec:symcost}. Powering is different: $D^{\mathrm{sym}}(a \odot p,a \odot q) = h(av)$ changes the scale at which the loss turns from quadratic to linear, with $h(av)/a^2$ tending to the quadratic CoDA regime as $a \to 0$ and $h(av)/a$ to the Hilbert range as $a \to \infty$; powering thus defines a parametric Huber family rather than acting as an equivariance.

Vertex data are excluded for a different reason than in case (b) but as clearly: $D^{\mathrm{sym}}(q,e_i) = +\infty$ for every interior $q$, through the term $\log\big(\frac{1}{n}\sum_jq_j/x_j\big)$, and the divergent term depends on $q$, so dropping data-side terms gives no finite objective: the device behind the categorical version \eqref{eq:lambda-mle} has no analogue here. Exact zeros are excluded; behaviour \emph{near} them is treated in \autoref{subsubsec:rob-sym}.

\autoref{tab:fits} collects the fits of the running example. Two things it makes visible: EKO and the two $\lambda$-estimators are fitted on one-parameter families, the Fr\'echet means and the barycenter on the saturated model, so the rows differ by model as well as by principle; and the flat Fr\'echet mean and the barycenter have no categorical row, both losses being infinite at a vertex.

\begin{table}[t]
\centering
\begin{tabular}{llll}
\hline
estimator & model & data & fit\\
\hline
EKO, $U = \exp$ & $E_t$ & either, through $\bar x$ & $(0.438,\,0.323,\,0.238)$\\
WZTY, $\lambda = 1$ & $E_t^{1}$ & interior sample & $(0.409,\,0.333,\,0.258)$\\
$\lambda$-ML, $\lambda = 1$ & $E_t^{1}$ & categorical counts & $(5/12,\,4/12,\,3/12)$\\
LCF, Fisher-Rao & saturated & interior sample & $(0.503,\,0.196,\,0.301)$\\
LCF, Fisher-Rao & saturated & categorical counts & $(0.599,\,0.135,\,0.266)$\\
LCF, flat & saturated & interior sample & $(0.506,\,0.191,\,0.303)$\\
$D^{\mathrm{sym}}$-barycenter & saturated & interior sample & $(0.504,\,0.192,\,0.304)$\\
\hline
\end{tabular}
\caption{The fits on the running example, $n = 3$ and $t(i) = i$: the interior sample $x_1,x_2,x_3$ of \autoref{subsec:eko} and the categorical counts with the same mean $\bar x = (0.5,0.2,0.3)$. $E_t$ is the exponential family, $E_t^{1}$ the $\lambda = 1$ family \eqref{eq:lambda1-family}, saturated the whole simplex; WZTY uses the uniform geometric-mean generator; for categorical counts the EKO row is maximum likelihood; the flat row is the CoDA center.}
\label{tab:fits}
\end{table}

\subsection{Robustness and efficiency}\label{subsec:robustness}

We now compare the three principles in terms of stability and precision, treating each estimator as an M-functional of the data law, $T(F) = \arg\min_{q\in M}\E_F[\ell(q,X)]$, with \eqref{eq:m-estimator} its empirical version \citep{huber1964,vandervaart1998}. Two conventions hold throughout. Comparisons are made at a fixed estimand, the working family and the quantity it parametrizes held fixed while the loss varies and at laws where the estimands coincide, as in the matched-versus-decoupled density-power comparison of \autoref{subsec:eko}; and influence is measured in the parameter norm for parametric estimators, in the metric norm for the Fr\'echet mean, and through a common functional when the principles are compared with each other.

\subsubsection{The perturbation framework}\label{subsubsec:rob-framework}

The influence function measures the effect on the estimate of an infinitesimal contamination at $x$,
\begin{equation}\label{eq:ef-def}
  \mathrm{IF}(x) := \lim_{\varepsilon\downarrow0}
  \frac{T\big((1 - \varepsilon)F + \varepsilon\delta_x\big) - T(F)}{\varepsilon};
\end{equation}
its supremum over $x$ is the gross-error sensitivity, and bounded influence means that no single observation can move the estimate by more than a fixed amount. Write $\psi(x,q) = \nabla_q\ell(q,x)$ and $H = \nabla_q\E_F[\psi(X,q)]$ for the risk Hessian at the fit. Standard M-estimation calculus gives the influence function and, by the delta method, the asymptotic covariance,
\begin{equation}\label{eq:if-sandwich}
\begin{aligned}
  \mathrm{IF}(x)&= -H^{-1}\,\psi\big(x,T(F)\big),
  \qquad
  \Omega = \Cov_F\big(\psi(X,T(F))\big),\\
  \mathrm{Avar}&= H^{-1}\,\Omega\,H^{-\top}
  = \int \mathrm{IF}\,\mathrm{IF}\T\,dF.
\end{aligned}
\end{equation}
The formulas apply directly to the additive losses.

For the Bregman class everything reduces to the mean. Write $\mu: \mathcal P(\simplex) \to \simplex$, $\mu(F) = \E_F[X]$, for the mean map, and $\Pi_M: \simplex \to M$, $m \mapsto \arg\min_{q\in M}D_\psi[m:q]$, for the divergence projection, single-valued on the $U$-models by \autoref{thm:pyth}.

\begin{proposition}[Factorization of Bregman estimation]\label{prop:bregman-factor}
For the loss $\ell(q,x) = D_\psi[x:q]$ of a differentiable convex $\psi$ and any model $M$ on which $\Pi_M$ is defined, $T(F) = \Pi_M(\mu(F))$ for every $F$. If $\Pi_M$ is differentiable at $\mu(F)$, with derivative $A$ on the tangent space, then
\begin{equation}\label{eq:ef-affin}
  \mathrm{IF}(x) = A\,\big(x - \mu(F)\big),
\end{equation}
affine in the contamination point.
\end{proposition}

The proposition separates estimators where the mean suffices from those that retain information from the sample. Throughout the Bregman case, the contamination enters affinely and all loss-dependence enters via the constant gain $A$; the $\lambda$-projection and the Fr\'echet estimator instead depend nonlinearly on the observation, through the score weights and the logarithm map respectively.

The gross-error sensitivity is $\gamma^\ast(T,F) = \sup_x\|\mathrm{IF}(x)\|$ in the norm just fixed. In the Bregman class affinity makes the supremum occur at a vertex; away from that case, the question is whether influence remains bounded as contamination approaches the boundary.

\subsubsection{The dually flat case}\label{subsubsec:rob-eko}

Fisher consistency holds throughout: each $D_U[\cdot:q]$ is a divergence in its second argument, so whenever the mean of the data law lies on the family, $\mu(F) = p_{\theta_0}$, the factorization proposition returns $T(F) = \Pi_M(p_{\theta_0}) = p_{\theta_0}$ for every admissible $U$, and all admissible $U$ estimate the same point.

For the matched pair $(D_U,E_t^U)$ the stationarity condition is the unweighted moment matching \eqref{eq:u-momentenabgleich}, so the factor map solves $\E_{p_\theta}[t] = \langle m,t\rangle$ in $\theta$ and $\mathrm{IF}_\theta(x) = G_U^{-1}(\langle x,t\rangle - \eta)$ with $G_U = \nabla_\theta\E_{p_\theta}[t]$, $\eta = \E_{p_{\theta_0}}[t]$: the observation enters through the functional $\langle x,t\rangle$ at every $\beta$, the weights cancelling as in \eqref{eq:beta-kuerzung}, and the $U$-dependence survives only in the gain $G_U^{-1}$. Robustness would have to show up as observation-dependent damping of the contamination, and on the matched pair changing $\beta$ introduces none: contamination enters through the same linear functional $\langle x,t\rangle$ at every $\beta$, while $\beta$ modifies only the gain $G_U^{-1}$. This is the quantified form of the claim of \autoref{subsec:eko}.

In the decoupled situation, robustness reappears. For the density-power loss on the undeformed family $E_t$, the sandwich formulas \eqref{eq:if-sandwich} evaluated at the model, with the $\beta$-tilted information $J_\beta = \sum_i p(i)^{1+\beta}\big(t(i) - \eta\big)\big(t(i) - \eta\big)\T$ as risk Hessian, gives (\autoref{app:dpd-score})
\begin{equation}\label{eq:ef-dpd}
  \mathrm{IF}_\beta(x)
  = J_\beta^{-1}\Big(\big\langle x,\ p^{\beta}(t - \eta)\big\rangle
  - \zeta_\beta\Big),
  \qquad \zeta_\beta = \sum_i p(i)^{1+\beta}\big(t(i) - \eta\big),
\end{equation}
still affine in the contamination point, as \autoref{prop:bregman-factor} requires. What $\beta$ changes is that the observation now enters through $\langle x,\,p^{\beta}(t - \eta)\rangle$ rather than through the raw $\langle x,\,t - \eta\rangle$ of the matched case: the $j$-th coefficient carries the damping $p(j)^{\beta}$, so that influence decreases with the cell probability, now with respect to the directions (the parts on which $x$ places its mass) as opposed to the locations (how far $x$ lies from the fit). Since $\mathrm{IF}$ is affine in $x$, its norm is convex and attains its maximum over the compact simplex at a vertex.

There is a trade-off with efficiency. At the model the score is centred and $\mathrm{Avar}_\beta = J_\beta^{-1}a_\beta\T\Cov_F(X)a_\beta J_\beta^{-\top}$ with $a_\beta = p^{\beta}(t - \eta)$, the sampling variability entering through the covariance of the composition law. The categorical case fixes that covariance at $\diag(p) - pp\T$, where $\beta = 0$ attains the Cramér-Rao bound of the multinomial \citep{basu1998} and every larger $\beta$ may incur an efficiency loss; interior data leave the covariance free, and how much efficiency a given $\beta$ costs then depends on the dispersion of the data. The running sample shows the expected robustness-efficiency trade-off, quantified in \autoref{app:rob-eko-numerics}. The $\gamma$-divergence of Fujisawa and Eguchi \citep{fujisawaeguchi2008} is a boundary case: although its logarithm sits outside the sample average, its risk still factors through $\bar x$.

\subsubsection{The projectively flat case}\label{subsubsec:rob-wzty}

WZTY does not keep one estimand while the generator varies: as \autoref{subsec:wzty} noted, the pair $(\phi,E_t^\lambda)$ defines the population projection $T(F) = \arg\min_\theta\E_F D^{(\lambda)}_\phi(X \mid p_\theta)$ itself, and for a fixed pair the perturbation theory is standard M-estimation. With $J_\theta$ and $q = p_\theta$ as in \autoref{subsec:wzty}, the per-observation score of \eqref{eq:lambda-score} is $\psi(x,\theta) = -\,J_\theta\T\,w(x,q)\,g(q)\,(x - q)$ with $w(x,q) = \big(1 + \lambda\,\nabla\phi(q) \cdot (x - q)\big)^{-1}$, and \eqref{eq:if-sandwich} gives $\mathrm{IF}(x) = -H^{-1}\psi(x,\theta_0)$. The shape is an affine profile $g(q)(x - q)$ modulated by the scalar weight $w$, hence a dependence on where the observation is located, which cannot occur in the \autoref{prop:bregman-factor} case. As $\lambda \to 0$, $w \to 1$ and $g \to -\nabla^2\phi$, recovering the mean-factoring class.

At $\lambda = 1$ with the weighted-geometric-mean generator the modulation has a closed form (using Aitchison operations \autoref{subsec:Aitch_def}). There $w(x,q) = \big(\sum_l\pi_l x_l/q_l\big)^{-1}$ and, componentwise,
\begin{equation}\label{eq:lambda-residuum}
  \big(-\,w\,g(q)(x - q)\big)_i
  = \frac{\pi_i - \hat r_i(x,q)}{q_i},
  \qquad
  \hat r_i(x,q) = \frac{\pi_i\,x_i/q_i}{\sum_{l=1}^n \pi_l\,x_l/q_l},
\end{equation}
the closure of the $\pi$-weighted likelihood ratios: $\hat r = \pi \oplus (x \ominus q)$, the Aitchison residual of the observation from the fit, then perturbed onto the reference weights. The influence of an observation is carried entirely by this residual composition; contamination at $x = q$ has $\hat r = \pi$ and vanishing score.

As $x \to e_j$ the weight stays bounded, $w \to q_j/\pi_j$, and the $j$-th entry of the score core $-\,w\,g(q)(x - q)$ grows like $1/q_j$ as the fitted $q_j \to 0$: rare parts of the fit are amplified, opposite to the $p^{\beta}$ damping of the decoupled density-power loss and consistent with the escort factor $p^{-\lambda}$ in \eqref{eq:begleit-moment}. The logarithmic projection responds to ratios whereby the Bregman class responds to differences.
The one-dimensional case shows the complementary bounded-tail mechanism explicitly. For $n = 2$ with uniform weights, $\phi(q) = \tfrac{1}{2}\log q_1 + \tfrac{1}{2}\log q_2$, in the log-odds coordinate $z = \log(q_1/q_2)$ that simplifies the flat geometry of \eqref{eq:ldiv-clr-rueckzug}, the $L$-divergence is
\begin{equation}\label{eq:logcosh}
  D_\phi(x \mid q) = \log\cosh\Big(\frac{z_x - z_q}{2}\Big),
\end{equation}
a log-cosh loss with quadratic core $\tfrac{1}{8}(z_x - z_q)^2 + O((z_x - z_q)^4)$, linear tails $\tfrac{1}{2}|z_x - z_q| - \log2$, and score $\tfrac{1}{2}\tanh((z_q - z_x)/2)$ (\autoref{app:rob-wzty-logcosh}). Its Huber crossover is fixed by the geometry rather than chosen, and what it regularizes is the squared Aitchison loss of case (b) of \autoref{subsec:lcf}, proportional to $(z_x - z_q)^2$ in this coordinate: that loss has score linear in $z_q - z_x$, so in the unbounded clr chart one remote observation moves the estimate without limit, while the log-cosh keeps the same quadratic behaviour near the diagonal and flattens the tail.

\begin{proposition}[Bounded response sensitivity of the logarithmic projection]
\label{prop:wzty-bounded}
Fix an interior weight $\pi$ and an interior fit $q$. The per-observation score of the $L$-divergence loss has free components $\psi_j(x,q) = (\pi_j/q_j)\big(1 - r_j/\sum_i\pi_ir_i\big)$, $r_i = x_i/q_i$, bounded over the closed simplex, the supremum of the bracket being $\max\big(1,\,1/\pi_j - 1\big)$, attained in the limits $x_j \to 0$ and $x \to e_j$. Hence the gross-error sensitivity is finite at every interior fit but its bound grows as $1/q_j$ toward rare parts: robustness is pointwise in the fitted parameter, not uniform.
\end{proposition}

\subsubsection{The metric case}\label{subsubsec:rob-lcf}

For LCF, $T(F) = \arg\min_q\E_F d(q,X)^2$ is itself a point of the simplex and influence is measured in the metric norm. By \autoref{prop:grad-dist}, $\psi(x,q) = -2\log_q(x)$ and population stationarity is $\E_F\log_{q_*}(X) = 0$. Differentiating this along the contaminated path gives \eqref{eq:if-sandwich} in geometric form,
\begin{equation}\label{eq:if-lcf}
  \mathrm{IF}(x) = 2\,H^{-1}\log_{q_{\ast}}(x),
\end{equation}
mean zero by stationarity, with $H$ the covariant Hessian of the Fr\'echet objective at its minimum, an operator on $T_{q_{\ast}}\simplex$; the asymptotic covariance is the bilinear form $H^{-1}\Omega H^{-1}$ on that tangent plane, $\Omega = 4\,\Cov_F\big(\log_{q_{\ast}}(X)\big)$. Constant curvature makes $H$ explicit, from the eigenvalues \eqref{eq:sphere-hessian-ev}: writing $d_k = d(q_{\ast},x_k)$ and $u_k = \log_{q_{\ast}}(x_k)/d_k$ for the unit direction toward the $k$-th observation,
\begin{equation}\label{eq:frechet-hesse}
  H = \frac{2}{N}\sum_{k=1}^N\Big(u_ku_k\T
  + c_k\big(I - u_ku_k\T\big)\Big),
  \qquad
  c_k = \sqrt{\kappa}\,d_k\cot(\sqrt{\kappa}\,d_k),
\end{equation}
the population version obtained by integrating over $F$. The comparison function that gave uniqueness in \autoref{subsec:lcf} is thus also the gain of the perturbation theory.

Two elements control the influence: the numerator has $g$-length exactly $d(q_{\ast},x)$, bounded by the diameter of \eqref{eq:diameter}; the gain $H^{-1}$ is bounded by $1/\lambda_{\min}(H)$, and every $c_k$ is positive when $d_k < \pi/(2\sqrt{\kappa})$, which for Fisher-Rao is the diameter $\pi$, so for interior data the empirical Hessian is positive definite and uniform control of its inverse is a lower bound on its smallest eigenvalue. Together,
\begin{equation}\label{eq:if-lcf-bound}
  \big\|\mathrm{IF}(x)\big\|_{g}
  \ \le\ \frac{2\,d(q_{\ast},x)}{\lambda_{\min}(H)}
  \ \le\ \frac{2\pi}{\lambda_{\min}(H)}:
\end{equation}
distance bounds the displacement and curvature determines the gain. For Fisher-Rao a rare fitted part pushes $d(q_*,e_j) = 2\arccos\sqrt{(q_*)_j}$ toward the diameter $\pi$ without reaching it (numerical coefficients and influence norms in \autoref{app:rob-lcf-numerics}). In the flat case, $c \equiv 1$, $H = 2I$ and $\mathrm{IF}(x) = \log_{q_*}(x)$, whose norm is the Aitchison distance: finite for every interior contamination and divergent at the boundary, so the zero problem reappears (as a perturbation) as infinite gross-error sensitivity. Nonconstant curvature is handled by the same distance-Hessian factorization through the Jacobi equation. The score $\log_q(x)$ is single-valued only where $x$ is joined to $q$ by a unique minimizing geodesic, and the injectivity radius is the distance up to which that holds: within that radius, the Hessian of squared distance is block diagonal in the splitting $T_p\simplex = \R\,\dot\gamma(r) \oplus \dot\gamma(r)^{\perp}$ at $r = d(p,x)$ \citep[Ch.~6]{jost2017}, \citep[\S1.5]{karcher1977},
\begin{equation}\label{eq:hessian-blocks}
\Hess_p\,\tfrac{1}{2} d(\cdot,x)^2
\;=\;1\,\oplus\,r\,J'(r)\,J(r)^{-1},
\end{equation}
with $J$ the matrix Jacobi field along the minimizing geodesic from $x$ to $p$ (\autoref{app:rob-geometric-preliminaries}); the normal block is defined below the first conjugate distance, where $J(r)$ is invertible,\footnote{A Jacobi field with $J(0) = 0$ vanishes again where the geodesics from $x$ refocus; the smallest such radius is the first conjugate distance, and $J(r)$ is singular there, so the inverse exists strictly below it. For Fisher-Rao it is $2\pi$, against a diameter of $\pi$.} and positive definite at least below the comparison radius $\pi/(2\sqrt{K_{\max}})$ when $K_{\max} > 0$, since $R \preceq K_{\max}I$ along $\gamma$ gives $r\,J'(r)J(r)^{-1} \succeq \sqrt{K_{\max}}\,r\cot(\sqrt{K_{\max}}\,r)\,I$. One transverse factor therefore controls both strict convexity and the influence gain, the constant-curvature coefficients $c_k$ of \eqref{eq:frechet-hesse} being its special case.

To wit, in the square-root geometry of \autoref{subsec:lcf} the curvature extends smoothly to the compactification of \autoref{prop:sqrt-sphere} and the injectivity radius is positive, so \eqref{eq:hessian-blocks} applies there, with $K_{\max} = 5/4$ for $n = 3$; past the comparison radius, positivity of the vertex Hessian is a computational task rather than a global uniqueness theorem (\autoref{app:rob-sqrt-certificates}).

\subsubsection{The symmetrized case}\label{subsubsec:rob-sym}

The barycenter of \autoref{subsec:jbary} is an additive M-estimator, so \eqref{eq:if-sandwich} applies as given, and it gives more concrete bounds than the diameter-led quantities above. Working in the clr chart, where the functional is $T(F) = \arg\min_\xi\E_F\,h\big(\xi - \clr(X)\big)$, the per-observation score is $\psi(x,\xi) = \nabla h\big(\xi - \clr(x)\big) = \clr^{-1}\big(\xi - \clr(x)\big) - \clr^{-1}\big(\clr(x) - \xi\big)$, the risk Hessian is $M = \E_F\,\nabla^{2}h\big(\xi_{\ast} - \clr(X)\big)$, positive definite on $V_0$ for every $F$ since $\nabla^{2}h(v)$ has kernel exactly $\operatorname{span}(\one)$ at every $v$, and
\begin{equation}\label{eq:ef-jbary}
  \mathrm{IF}(x)
  = -\,M^{-1}\Big(\clr^{-1}\big(\xi_{\ast} - \clr(x)\big)
  - \clr^{-1}\big(\clr(x) - \xi_{\ast}\big)\Big),
  \qquad
  \big\|\mathrm{IF}(x)\big\| \le \frac{2}{\lambda_{\min}(M)},
\end{equation}
the bound holding over the closed simplex: the score is a difference of two probability vectors wherever it is defined, and its limits as $x$ degenerates towards the boundary are of the same form. The gross-error sensitivity is finite.

A comparison with the flat LCF case makes the mechanism patent: the induced metric of $D^{\mathrm{sym}}$ is the same flat Aitchison-type metric up to scale, yet its Fr\'echet mean has unbounded clr influence while \eqref{eq:ef-jbary} is bounded. What matters is a bounded loss gradient: finite diameter produces one for squared distance, while the Hilbert-slope tails \eqref{eq:hilbert-sandwich} produce one without curvature or compactness. Finite diameter is therefore sufficient, not necessary. The separation is visible under contamination $x_4 = C(1,\delta,\delta)$ drifting to a vertex: the $D^{\mathrm{sym}}$-barycenter converges to the finite limit $\approx (0.702,0.122,0.177)$ while the geometric mean diverges logarithmically (\autoref{app:rob-sym-numerics}). Stability near the boundary and existence on it are thus different properties: the barycenter has bounded gross-error sensitivity although the loss is infinite at the vertex itself. By the quadratic part \eqref{eq:symcost-germ} the influence function approaches the CoDA center's when both the data and the contamination remain in the local quadratic regime, and the Huber-type trade-off opens only with dispersion, when observations reach the linear tails and the score settles at its bound.

\subsubsection{Synthesis}\label{subsubsec:rob-synthesis}

Efficiency comparisons require a common estimand, and robustness requires a common norm. If $m(q) = \E_q[t]$, then $\mathrm{IF}_{m\circ T} = dm(\mathrm{IF}_T)$; an estimator bounded in its intrinsic metric can therefore be amplified by a functional emphasizing rare parts, and since EKO, WZTY and LCF in general estimate different functionals, a numerical ranking without this normalization is not meaningful. The mechanisms themselves are distinct: the arithmetic mean suffices for Bregman estimation and thus it has affine influence; WZTY modulates each observation by the likelihood-ratio weight $w_k$; LCF has score norm $2d(q,x)$, with diameter bounding the displacement and sectional curvature controlling the Hessian gain; and the symmetrized loss shows that curvature is not required for bounded influence, its flat metric being unbounded while its linear tails bound the loss gradient. In particular the logarithmic estimator is robust for its own $L$-projection functional, not a robust recovery of the arithmetic mean.

\section{Applications: estimation as a choice of geometry}\label{sec:Applications}

The geometries assembled above are the structures that applied work on categorical data already operates with, mostly without naming them. An estimation procedure on the simplex is a geometric choice, a metric together with an affine connection, or a divergence standing in for both, \textbf{and different choices return different results from the same data}; the choice is made whether or not it is declared. We make this explicit for difference-in-differences on categorical outcomes, where building a counterfactual composition from three observed ones is the choice of a parallelism on the simplex and quantify, on a designed example, how far the reported effect moves with it. A shorter second instance, the entropy-minimizing pricing measure of mathematical finance read as an exponential projection in the sense of \autoref{subsec:eko}, is in the Supplement (\autoref{app:follmer}).

Categorical data arrive as positive quantities $q \in \R^n_{>0}$; as in \autoref{sec:Intro}, the simplex is the quotient of this cone by the choice of num\'eraire, and a categorical treatment effect splits accordingly into an aggregate part, the effect on the total, a scalar problem for which log-linear difference-in-differences is adequate, and a compositional part on the base. Everything below acts on the base; if total effects are needed, an assumption on the fiber\footnote{Closure \eqref{eq:closure} presents the cone as a bundle over the simplex: the base is the set of compositions, and the fiber over one of them is the ray of its positive multiples, coordinatized by the total $\sum_i q_i$. A compositional effect moves along the base, an effect on the total along this fiber.} is needed, for instance the parallel-growths assumption of \citep{boussim2025}.

\subsection{Compositional difference-in-differences and the choice of transport}\label{subsec:app-cid}

Difference-in-differences on categorical outcomes provides four compositions: the control group before and after, $\omega_{00},\omega_{01} \in \simplex$\footnote{We write the cells as $\omega$; Boussim \citep{boussim2025} writes $\pi$, which here is reserved for the generator weights of \autoref{subsec:generator}.}, the treated group before, $\omega_{10}$, and the treated group after, $\omega^{I}_{11}$. The effect of the treatment is the comparison of $\omega^{I}_{11}$ with the never observed $\omega^{N}_{11}$, the composition the treated group would have shown untreated. Every procedure of this kind builds that fourth composition from the three observed ones by one rule: the treated group moves parallel to the control group.

Parallel in which sense is what the rule leaves open, since the change is recorded at one composition and applied at another. In $\R^n$ the question does not arise: the control change $\omega_{01} - \omega_{00}$ is a vector, and the same vector is available at every point. The simplex has that translation available too, as an open subset of an affine hyperplane, and it is one parallelism among several, the mixture connection, whose translate can leave the simplex; the log-ratio calculus supplies another. Which of them carries the control change from $T_{\omega_{00}}\simplex$ to $T_{\omega_{10}}\simplex$ has to be declared before the rule says anything, and an affine connection is what declares it. The identifying assumption of the design is thus not parallel trends by themselves, but parallel trends for a declared parallelism.

Boussim \citep{boussim2025} supplies one by working in the perturbation calculus of \autoref{subsec:Aitch_def}. The control change is read off as the perturbation difference $\omega_{01} \ominus \omega_{00} = C(\omega_{01}/\omega_{00})$ and applied at the treated baseline by perturbation:
\begin{equation}\label{eq:boussim-cf}
\omega^{N}_{11} = \omega_{10} \oplus (\omega_{01} \ominus \omega_{00}) = C\!\left(\frac{\omega_{10}\,\omega_{01}}{\omega_{00}}\right).
\end{equation}
On log-ratios the operation is additive, and for any two parts $i,j$ it reads
\begin{equation}\label{eq:boussim-logquoten}
\log\frac{(\omega^{N}_{11})_i}{(\omega^{N}_{11})_j}
= \log\frac{(\omega_{10})_i}{(\omega_{10})_j}
+ \left(\log\frac{(\omega_{01})_i}{(\omega_{01})_j} - \log\frac{(\omega_{00})_i}{(\omega_{00})_j}\right),
\end{equation}
the ordinary two-way difference performed on each log-odds separately, with closure restoring normalization at the end. It is a repair of the raw-proportions recipe: the right-hand side of \eqref{eq:boussim-cf} lies in the open simplex for any four cells, while $\omega_{10} + (\omega_{01} - \omega_{00})$ need not.

Two readings of \eqref{eq:boussim-cf} are available, and both are in use. In the first the parallelism belongs to the Aitchison metric: by \eqref{eq:clr-isometrie} the map $\clr: \simplex \to V_0$ is an isometry onto $V_0 = \{v \in \R^n: \sum_iv_i = 0\}$ with its Euclidean structure, so Levi-Civita transport is translation in that chart, and since $\clr(p \oplus q) = \clr(p) + \clr(q)$ the perturbation in \eqref{eq:boussim-cf} is that translation. In the second the parallelism is the exponential connection of the Fisher-Rao dually flat pair, flat in the coordinates $\operatorname{alr}$ of \eqref{eq:alr}, where transport is again translation and leaves every log-odds unchanged; that is \eqref{eq:boussim-logquoten}.

These are not two connections that happen to agree on one formula. Writing $y_i = \log p_i$,
\begin{equation}\label{eq:clr-alr-linear}
\clr(p) = \left(I - \tfrac{1}{n}\one\one\T\right)y,
\qquad
\operatorname{alr}(p) = (y_i - y_n)_{i<n},
\end{equation}
both linear in $y$ and both killing $\one$, so the transition map $L = \operatorname{alr} \circ \clr^{-1}: V_0 \to \R^{n-1}$, $L(v) = (v_i - v_n)_{i<n}$, is a linear isomorphism: $L(v) = 0$ forces $v_1 = \dots = v_n$, hence $v = 0$ on $V_0$. A connection flat in one chart is flat in every linearly related one, so the Levi-Civita connection of the Aitchison metric and the exponential connection of Fisher-Rao are one flat connection on the open simplex, seen in two log-ratio charts differing only in the num\'eraire, geometric mean against a single part, as noted at \eqref{eq:alr}.

What \eqref{eq:boussim-cf} does not fix is the metric, and a geometry in the sense of \autoref{subsec:generator} is a metric together with a connection. That one flat connection is the Levi-Civita connection of $g^A$, where both Eguchi connections vanish in the clr chart (\autoref{subsubsec:symcost}), and it is one half of the dual pair of $g^{FR}$, whose Levi-Civita connection is not it but the average \eqref{eq:lc-average}. The two metrics disagree about the simplex itself: Aitchison is flat with infinite diameter, Fisher-Rao has curvature $1/4$ and the finite diameter \eqref{eq:diameter}.

Building $\omega^{N}_{11}$ is half of the design; the other half is the comparison with $\omega^{I}_{11}$, and effect size, distance, the variance of an estimated cell and the mean of several cells are metric quantities, none of them fixed by \eqref{eq:boussim-cf}. Two analysts can agree on every counterfactual, cell by cell, and still report different effect sizes from the same four cells.

The procedure Boussim argues against completes the picture: linear parallel trends on raw proportions carry the control change by translation in the mixture coordinates, transport for the mixture connection of \autoref{subsec:generator}, and his Figure~1 is a design where the result is not a composition at all. The two procedures in the applied discussion are therefore the two halves of one dual pair, $\nabla^{(m)}$ against $\nabla^{(e)}$, and the discussion has run as one between a wrong answer and a right one. The pair has a third transport belonging to neither half, the Levi-Civita connection of $g^{FR}$, their average \eqref{eq:lc-average} and the parallelism of the metric under which Zhou et al.~\citep{zhou2025} treat the same problem; it is not flat, and it does not appear in the discussion.

Every construction of this kind is an instance of one formula,
\begin{equation}\label{eq:cid-template}
\omega^{\mathrm{cf}}
\;=\;
\exp^{\nabla}_{\omega_{10}}\!\left(\;\tau^{\nabla}_{\omega_{00}\to\omega_{10}}\,\log^{\nabla}_{\omega_{00}}(\omega_{01})\;\right),
\end{equation}
where $\log^{\nabla}_{\omega_{00}}: \simplex \to T_{\omega_{00}}\simplex$ writes the observed control change as a tangent vector, $\tau^{\nabla}_{\omega_{00}\to\omega_{10}}: T_{\omega_{00}}\simplex \to T_{\omega_{10}}\simplex$ is the identification called for above, and $\exp^{\nabla}_{\omega_{10}}$ returns it to the simplex. The free parameter is $\nabla$, an affine connection, plus a metric when $\nabla$ is Levi-Civita, plus a path when $\nabla$ is curved. The identification assumption of the design is exactly that the treated group's untreated change is the parallel-transported control change for the declared $\nabla$; changing $\nabla$ changes the assumption, not the computation.

Four instantiations besides the raw-proportions transport span the cases of \autoref{sec:Estimators}. First, the e-connection, where \eqref{eq:cid-template} returns \eqref{eq:boussim-cf} (his Theorem~1) and parallel trends hold in the log-odds \eqref{eq:boussim-logquoten} (his Proposition~1); the connection is flat, so no path enters. Second, the Levi-Civita connection of Fisher-Rao, carried out on the sphere of radius $2$ through the isometry \eqref{eq:sphere-iso}, where $\log$, $\tau$ and $\exp$ are elementary (\autoref{app:cid-curvature}); curvature $\kappa = \tfrac{1}{4}$ makes \eqref{eq:cid-template} path-dependent, and we run both the direct route $\omega_{00} \to \omega_{10}$ and the detour through $\omega_{01}$. Third, the projection transport of Zhou et al.~\citep{zhou2025} on the same sphere, which projects the control direction at $\sqrt{\omega_{00}}$ onto the tangent space at $\sqrt{\omega_{10}}$ and restores its length. It is not Levi-Civita transport: before restoring length, projection shortens the component along the geodesic joining the base points by the cosine of the angle between them and leaves the orthogonal component alone, and the rescaling then changes both, so the two agree only in degenerate configurations (\autoref{app:cid-curvature}), and its scalar instance, the change-in-changes push-forward, does not preserve length at all. Declaring the metric is thus not yet declaring the estimator. Fourth, the Levi-Civita connection of the square-root-generator metric \eqref{eq:wurzel-metrik}, whose curvature ranges over $[-\tfrac{1}{4},\tfrac{5}{4}]$ for $n = 3$; no closed form exists, and \eqref{eq:cid-template} is integrated in the fourth-root chart, checked against Gauss-Bonnet.

\begin{table}[t]
\centering
\begin{tabular}{lrrrrrr}
\hline
 & \multicolumn{3}{c}{$\omega^{\mathrm{cf}}$} & \multicolumn{3}{c}{CTT $\;=\;\omega^{I}_{11} \ominus \omega^{\mathrm{cf}}$}\\
transport ($\nabla$) & \multicolumn{1}{c}{$1$} & \multicolumn{1}{c}{$2$} & \multicolumn{1}{c}{$3$} & \multicolumn{1}{c}{$1$} & \multicolumn{1}{c}{$2$} & \multicolumn{1}{c}{$3$}\\
\hline
m-connection (raw proportions) & $-0.2000$ & $0.4000$ & $0.8000$ & & & \\
e-connection \citep{boussim2025} & $0.0338$ & $0.1775$ & $0.7887$ & $0.6300$ & $0.2800$ & $0.0900$\\
Fisher-Rao Levi-Civita & $0.0012$ & $0.2640$ & $0.7348$ & $0.9839$ & $0.0107$ & $0.0055$\\
projection transport \citep{zhou2025} & $0.0012$ & $0.2654$ & $0.7334$ & $0.9843$ & $0.0104$ & $0.0054$\\
square-root Levi-Civita & $0.0231$ & $0.2345$ & $0.7424$ & $0.7498$ & $0.1724$ & $0.0778$\\
\hline
\end{tabular}
\caption{Five parallelisms, one design. Cells from Boussim's Figure~1 \citep{boussim2025}: $\omega_{00} = (0.7,0.2,0.1)$, $\omega_{01} = (0.3,0.3,0.4)$, $\omega_{10} = (0.2,0.3,0.5)$, with $\omega^{I}_{11} = (0.15,0.35,0.50)$. The m-connection counterfactual leaves the simplex, so no compositional effect can be formed. Curved transports are along the direct route $\omega_{00} \to \omega_{10}$.}
\label{tab:cid-parallelisms}
\end{table}

\autoref{tab:cid-parallelisms} runs all five on Boussim's own illustration. Each row builds its counterfactual with its own $\nabla$, while the effect column applies one and the same difference to all of them, the perturbation $\omega^{I}_{11} \ominus \omega^{\mathrm{cf}}$, so that the rows can be read against each other as compositions. Each row could equally use the difference of its own geometry, for the curved ones the tangent vector $\log^{\nabla}_{\omega^{\mathrm{cf}}}(\omega^{I}_{11})$, one more place where the declared parallelism enters after the counterfactual has been built. The reported effect is not robust to that parallelism: the counterfactual first share moves by a factor of about $28$ between the e-connection ($0.0338$) and the Fisher-Rao Levi-Civita answer ($0.0012$), and the effect on the first category, $0.630$ for the e-connection, $0.984$ for the two sphere transports and $0.750$ for the square-root metric, is not one finding in three variants. The Levi-Civita and projection counterfactuals agree here to four decimals in the first share and differ by about $0.0014$ in the other shares, and path dependence is small here, both for one reason: the three base cells are within $4 \times 10^{-4}$ of geodesic collinearity, so the Fisher-Rao loop rotation is only $0.150^\circ$ and the direction gap between the two transports $0.227^\circ$, against about $4^\circ$ for both on the curved design of \autoref{tab:cid-curvature}. The design thus primarily illustrates the choice of connection, with little sensitivity to the transport path.

Curvature enters through holonomy. In a flat geometry transport around a closed loop returns every vector unchanged, so carrying the control change over is one instruction with one outcome; in a curved one, for the two-dimensional geometries considered here, the loop returns it rotated by $\int K\,dA$ over the enclosed surface, and the direct route $\omega_{00} \to \omega_{10}$ and the detour through $\omega_{01}$ deliver two counterfactuals from the same four cells and the same words. A curved-geometry difference-in-differences therefore carries a modelling choice the flat designs never had to mention, the path, and the holonomy angle measures it exactly: on cells spanning a proper geodesic triangle the loop rotation is about $4^\circ$ in both curved metrics and the direct and detour counterfactuals differ by about a fifth of the gap between the e- and Levi-Civita answers on those cells (\autoref{tab:cid-curvature}). Near the boundary the Fisher-Rao exponential map crosses the wall $\omega_1 = 0$ and the squared chart folds the geodesic back inside, so the returned composition is an artifact of the chart, while the square-root transport stays interior on the same designs (\autoref{app:cid-curvature}); a report should say on which side of this its data sit.

The statistical layer transfers without new work. For every $\nabla$ above the counterfactual is a smooth map $\Phi^{\nabla}: (\simplex)^{3} \to \simplex$ of the three observed cells, so influence functions compose by the chain rule and the analysis of \autoref{subsec:robustness} reaches each comparator through its own case, \autoref{subsubsec:rob-eko} for the e-connection and the Levi-Civita-Fr\'echet theory of \autoref{subsec:lcf} for the two Levi-Civita rows; the projection row likewise, though no principle of \autoref{sec:Estimators} produces it. When the cells are estimated, Zhou et al.~\citep{zhou2025} give consistency and rates by empirical-process arguments, with distributional inference open by their own account; when the cells are themselves Fr\'echet means of samples, their asymptotic covariance is the sandwich \eqref{eq:if-sandwich} with the covariant Hessian of \autoref{subsubsec:rob-lcf}, the central limit theorem for intrinsic means being that of Bhattacharya and Patrangenaru \citep{bhattacharyapatrangenaru2005}, under their regularity assumptions; the Fr\'echet variance of a cell, the minimum of its objective, has the limit of Dubey and M\"uller \citep{dubeymueller2019}.

\subsection{Outlook}\label{subsec:outlook}

Three threads continue elsewhere. The symmetrized transport cost and its barycenter leave open the efficiency curve and the breakdown point recorded in \autoref{subsubsec:rob-sym}. The characterization of which losses are the canonical divergences of their own geometries, through Osband's principle on the simplex, is developed in \citep{navasmest}. And the $\lambda$-divergences have a sequential side: writing each factor $A_k(q)^{-1/\lambda}$ of the $\lambda$-Gibbs posterior \eqref{eq:lambda-gibbs} as the Laplace transform of a Gamma-distributed temperature makes the score weights $w_k$ posterior means of those temperatures and the curvature $-\lambda$ their variance, which connects the projectively flat principle to sequential inference, where the likelihood ratio acts as an e-value; that connection is left to future work.

\section*{Funding and competing interests}
The author received no specific funding for this work and declares no competing interests.

\appendix
\section{Derivations for the affine structures}\label{app:affine-derivations}

The calculations here are those for the affine structures of \autoref{sec:Affine}, most of them for the $L$-divergence, and for these we lean on Pal and Wong \citep{palwong2016,palwong2018}. They also prove that the geometry of the $L$-divergence has constant negative curvature; the derivation in \autoref{app:ldiv-projective} writes out the Christoffel symbols and the full curvature tensor.

\subsection{Elementary calculations for the $L$-divergence}\label{app:ldiv-calculations}
Let $\Phi = e^\phi$ be concave. The supporting-hyperplane inequality gives
\[
\Phi(p) \le \Phi(q) + \nabla\Phi(q) \cdot (p - q).
\]
Since $\nabla\Phi = \Phi\nabla\phi$, this becomes
\[
e^{\phi(p)} \le e^{\phi(q)}\left(1 + \nabla\phi(q) \cdot (p - q)\right).
\]
Taking logarithms yields
\[
\phi(p) - \phi(q) \le \log\left(1 + \nabla\phi(q) \cdot (p - q)\right),
\]
so $D_\phi[p:q] \ge 0$, with equality only at $p = q$ when $\Phi$ is strictly concave.

For the weighted geometric mean $\phi(p) = \sum_i\pi_i\log p_i$, one has $\partial_i\phi(p) = \pi_i/p_i$ and
\[
1 + \nabla\phi(q) \cdot (p - q)
= 1 + \sum_i\frac{\pi_i}{q_i}(p_i - q_i)
= \sum_i\pi_i\frac{p_i}{q_i},
\]
using $\sum_i\pi_i = 1$. Substitution in \eqref{eq:l-divergence} gives \eqref{eq:ldiv-amgm}.

\subsection{Eguchi differentiation of the $L$-divergence}\label{app:ldiv-eguchi}
Write $A(p,q) = 1 + \nabla\phi(q) \cdot (p - q)$, so that
$D_\phi[p:q] = \log A(p,q) - \phi(p) + \phi(q)$. Differentiating once in the first argument,
\[
\partial_iD_\phi[p:q] = \frac{\partial_i\phi(q)}{A(p,q)} - \partial_i\phi(p).
\]
A derivative in the second argument gives
\[
\partial_{j'}\partial_iD_\phi[p:q]
= \frac{\partial_i\partial_j\phi(q)A - \partial_i\phi(q)\partial_{j'}A}{A^2},
\qquad
\partial_{j'}A = \sum_k\partial_j\partial_k\phi(q)(p_k - q_k) - \partial_j\phi(q).
\]
On the diagonal $q = p$, $A = 1$ and $\partial_{j'}A = -\partial_j\phi(p)$, hence
\[
\left.\partial_{j'}\partial_iD_\phi\right|_{q=p}
= \partial_i\partial_j\phi(p) + \partial_i\phi(p)\partial_j\phi(p),
\]
which, with the sign convention \eqref{eq:eguchi-metric}, is \eqref{eq:ldiv-metric}.

Differentiating once more in the first argument gives
\[
\partial_i\partial_jD_\phi[p:q]
= -\frac{\partial_i\phi(q)\partial_j\phi(q)}{A(p,q)^2} - \partial_i\partial_j\phi(p).
\]
Unlike its Bregman counterpart, this still depends on $q$ through $A$ and $\nabla\phi(q)$, so its $q$-derivative on the diagonal, the mixture Christoffel symbol in the coordinates $p$, is in general nonzero; \autoref{app:ldiv-projective} computes it, together with the curvature $-1$ of the connection, which is therefore not flat in any coordinates.

\subsection{Weighted clr pullback}\label{app:weighted-clr}
The map $\clr_\pi$ is a diffeomorphism of $\simplex$ onto $V_\pi$: it takes values in $V_\pi$ because $\sum_i\pi_i = 1$, and $z \mapsto C(e^{z_1},\dots,e^{z_n})$ inverts it on $V_\pi$, as for $\clr$ in \autoref{subsec:Aitch_def}. For $u \in T_p\simplex$ write $w_i = u_i/p_i$. Since $d(\log p_i)(u) = w_i$, the differential of \eqref{eq:clr-pi} removes the $\pi$-weighted mean,
\[
(d\clr_\pi(u))_i = w_i - \sum_j\pi_jw_j.
\]
Its squared norm in $V_\pi$ is therefore the $\pi$-weighted variance
\[
\begin{aligned}
\left\langle d\clr_\pi(u),d\clr_\pi(u)\right\rangle_\pi
&= \sum_i\pi_iw_i^2 - \left(\sum_i\pi_iw_i\right)^2\\
&= \sum_i\frac{\pi_i}{p_i^2}u_i^2 - \left(\sum_i\frac{\pi_i}{p_i}u_i\right)^2,
\end{aligned}
\]
which is exactly the quadratic form associated with \eqref{eq:ldiv-beispiel-metrik}. Hence \eqref{eq:ldiv-clr-rueckzug} follows.

\section{Projective flatness of the logarithmic-divergence
connection}\label{app:ldiv-projective}

We verify \eqref{eq:christoffel-proj} and \eqref{eq:ldiv-curvature} of
\autoref{subsec:wzty} directly from the $\lambda$-logarithmic divergence
\eqref{eq:lambda-ldiv}, for a general $\lambda$-exponentially concave
generator $\phi$, general $\lambda > 0$ and general $n$; the $L$-divergence
case is $\lambda = 1$ and needs no separate argument. Throughout, $\phi_i$
and $\phi_{ij}$ abbreviate the partial derivatives of $\phi$ at $p$, and
$a := -\lambda\,\nabla\phi(p)$.

Write $A(p,q) = 1 + \lambda\,\nabla\phi(q) \cdot (p - q)$, so that
$D^{(\lambda)}_\phi(p \mid q) = \frac{1}{\lambda}\log A - \phi(p) + \phi(q)$ and
$\partial_iA = \lambda\,\partial_i\phi(q)$. Differentiating twice in the
first argument gives, at general $(p,q)$,
\[
  \partial_i\partial_j D^{(\lambda)}_\phi(p \mid q)
  = -\lambda\,\frac{\partial_i\phi(q)\,\partial_j\phi(q)}{A(p,q)^2}
  - \partial_i\partial_j\phi(p).
\]
One further derivative in $q_k$, using $A(p,p) = 1$ and
$\partial_{k'}A|_{q=p} = -\lambda\,\partial_k\phi(p)$, gives at $q = p$
\[
  \partial_{k'}\partial_i\partial_j D^{(\lambda)}_\phi\big|_{q=p}
  = -\lambda\big(\phi_{ik}\phi_j + \phi_i\phi_{jk}
  + 2\lambda\,\phi_i\phi_j\phi_k\big),
\]
hence, by the sign convention of \eqref{eq:eguchi-m-connection},
\[
  \Gamma^{(m)}_{ij,k}
  = \lambda\big(\phi_{ik}\phi_j + \phi_i\phi_{jk}\big)
  + 2\lambda^{2}\,\phi_i\phi_j\phi_k.
\]
The metric \eqref{eq:lambda-metric} gives
$\phi_{ij} = -g_{ij} - \lambda\,\phi_i\phi_j$; substituting for $\phi_{ik}$
and $\phi_{jk}$, the two cubic terms so produced cancel $2\lambda^{2}\phi_i\phi_j\phi_k$ and, with
$a = -\lambda\nabla\phi$,
\[
  \Gamma^{(m)}_{ij,k} = a_i\,g_{jk} + a_j\,g_{ik},
\]
the projective form \eqref{eq:christoffel-proj}. Raising the index,
$\Gamma^{(m)\,l}_{ij} = \delta^l_i a_j + \delta^l_j a_i$.

For the curvature, since $a$ is a gradient,
$\partial_i a_j = \partial_j a_i$, and the metric relation above gives
$\partial_i a_k = -\lambda\,\phi_{ik} = \lambda\,g_{ik} + a_i a_k$.
Substituting
$\Gamma^{(m)\,l}_{ij} = \delta^l_i a_j + \delta^l_j a_i$ into the curvature
tensor
\[
  R^{(m)\,l}{}_{kij} = \partial_i\Gamma^{(m)\,l}_{jk}
  - \partial_j\Gamma^{(m)\,l}_{ik}
  + \Gamma^{(m)\,l}_{im}\Gamma^{(m)\,m}_{jk}
  - \Gamma^{(m)\,l}_{jm}\Gamma^{(m)\,m}_{ik}
\]
and using the symmetry of $\partial a$ to cancel the terms proportional
to $\delta^l_k$, the quadratic terms in $a$ combine to $\delta^l_i a_ja_k - \delta^l_j a_ia_k$ and
\[
  R^{(m)\,l}{}_{kij}
  = \delta^l_j\big(\partial_ia_k - a_ia_k\big)
  - \delta^l_i\big(\partial_ja_k - a_ja_k\big)
  = \lambda\,\big(\delta^l_j\,g_{ik} - \delta^l_i\,g_{jk}\big),
\]
the constant-curvature form
$R^{l}{}_{kij} = \kappa(\delta^l_ig_{jk} - \delta^l_jg_{ik})$ with
$\kappa = -\lambda$, the constant curvature found by Wong \citep{wong2018}; at
$\lambda = 1$ this is the value Pal and Wong \citep{palwong2018} prove for the $L$-divergence.

\section{Proofs for the symmetrized cost}\label{app:symmetry}

Proofs for the symmetrized cost of \autoref{subsubsec:symcost}, beginning with the inversion identity.

\begin{proof}[Proof of the inversion identity {$D_\phi[\iota(p):\iota(q)] = D_\phi[q:p]$}]
With $\Sigma(p) = \sum_ip_i^{-1}$ the ratios transform as
$\iota(p)_i/\iota(q)_i = (q_i/p_i)\,\Sigma(q)/\Sigma(p)$, so
\[
\begin{aligned}
\log\left(\frac{1}{n}\sum_i\frac{\iota(p)_i}{\iota(q)_i}\right)
&= \log\left(\frac{1}{n}\sum_i\frac{q_i}{p_i}\right) + \log\Sigma(q) - \log\Sigma(p),\\
\frac{1}{n}\sum_i\log\frac{\iota(p)_i}{\iota(q)_i}
&= \frac{1}{n}\sum_i\log\frac{q_i}{p_i} + \log\Sigma(q) - \log\Sigma(p).
\end{aligned}
\]
The $\Sigma$-terms cancel in the difference, leaving $D_\phi[q:p]$.
\end{proof}

\begin{proof}[Proof of \autoref{lem:hilbert-sandwich}]
From $\max_iv_i - \log n \le \Lambda(v) \le \max_iv_i$ applied to $v$ and $-v$,
$h(v)$ lies between the range of $v$ minus $2\log n$ and the range of $v$.
The sharper upper bound $2\log\cosh(d_{\mathrm{Hil}}/2)$ follows from the Kantorovich inequality
\citep{kantorovich1948}: for weights summing to one and ratios
$m \le r_i \le M$,
$\big(\sum_i\pi_ir_i\big)\big(\sum_i\pi_i/r_i\big) \le (M + m)^2/(4Mm)$, and
with $r_i = e^{v_i}$, $M/m = e^{d_{\mathrm{Hil}}}$, the right-hand side is
$\cosh^2(d_{\mathrm{Hil}}/2)$.
\end{proof}

\begin{proof}[Proof of \autoref{prop:symcost-degenerate}]
Compute in the clr chart, where $D^{\mathrm{sym}} = h(\xi - \eta)$ by
\eqref{eq:symkosten-clr}. The metric \eqref{eq:eguchi-metric} is
$-\partial_{\xi_i}\partial_{\eta_j}h(\xi - \eta)|_{\xi=\eta}
= \nabla^2h(0)_{ij}$, and \eqref{eq:symcost-hessian} at $v = 0$ gives
$\nabla^2h(0) = \frac{2}{n}\big(I - \frac{1}{n}\one\one\T\big)$, acting
on $V_0$ as $\frac{2}{n}$ times the identity; pulled back through clr this is
$\frac{2}{n}\,g^A$ by \eqref{eq:clr-isometrie}, twice the uniform-weight case of
\eqref{eq:ldiv-clr-rueckzug}. For the connections,
$-\partial_{\xi_i}\partial_{\xi_j}\partial_{\eta_k}
h(\xi - \eta)|_{\xi=\eta} = +\,h_{ijk}(0)$ and
$-\partial_{\eta_i}\partial_{\eta_j}\partial_{\xi_k}
h(\xi - \eta)|_{\xi=\eta} = -\,h_{ijk}(0)$, and every third derivative of
the even function $h$ vanishes at the origin.
\end{proof}

\begin{proof}[Proof of \autoref{prop:symcost-exclusion}]
Suppose $D^{\mathrm{sym}} = D_\psi$ in a chart $\theta$. Symmetry of the
divergence reads
$G(x,y) := 2\psi(x) - 2\psi(y) - \big(\nabla\psi(x) + \nabla\psi(y)\big) \cdot (x - y)
\equiv 0$ in the chart, and
$\partial_{x_k}\partial_{x_i}G = -\sum_j\psi_{ijk}(x)\,(x_j - y_j)$, which can
vanish for all $y$ only if $\nabla^3\psi \equiv 0$: the potential is
quadratic,
$D^{\mathrm{sym}}[p:q] = \tfrac{1}{2}(\theta_p - \theta_q)\T A\,(\theta_p - \theta_q)$
with $A = \nabla^2\psi$. The Eguchi metric of such a divergence is the
constant $A$ in the $\theta$-chart, with geodesic distance
$d(p,q)^2 = (\theta_p - \theta_q)\T A(\theta_p - \theta_q)$, so the
divergence equals half its own squared induced distance. The induced
metric is \eqref{eq:symcost-metric}, whose distance is
$d(p,q)^2 = \frac{2}{n}\|\xi_p - \xi_q\|^2$ by \eqref{eq:clr-isometrie};
hence $h(v) = \frac{1}{n}\|v\|^2$ identically, which contradicts
\eqref{eq:hilbert-sandwich}: its upper bound is linear along clr rays.
For the second claim, the mixture connection of $D^{(\lambda)}_\phi$ has,
by \autoref{prop:proj-flat}, curvature tensor of constant-curvature form
with $\kappa = -\lambda \ne 0$, while the mixture connection of
$D^{\mathrm{sym}}$ vanishes in the clr chart and its curvature is zero.
\end{proof}

\subsection{The Hessian in exponential-weight form}\label{app:hessian-weights}
The inverse clr map of \autoref{subsec:Aitch_def} sends $v \in V_0$ to $C(e^{v_1},\dots,e^{v_n})$, with $i$-th component $e^{v_i}/\sum_je^{v_j}$. With $\Lambda(v) = \log\big(\frac{1}{n}\sum_ie^{v_i}\big)$ as in \eqref{eq:symkosten-clr}, one differentiation gives $\nabla\Lambda(v) = \clr^{-1}(v)$ and a second gives $\nabla^2\Lambda(v) = \diag(q) - qq\T$ at $q = \clr^{-1}(v)$, the covariance matrix of the categorical law $q$. Since $h(v) = \Lambda(v) + \Lambda(-v)$ and $\clr^{-1}(-v) = \iota(\clr^{-1}(v))$ by \eqref{eq:inversion}, this is \eqref{eq:symcost-hessian}. In the weights $\sigma^{\pm} = \clr^{-1}(\pm v)$ it reads
\[
\nabla^2h(v) = \diag(\sigma^+) - \sigma^+(\sigma^+)\T
+ \diag(\sigma^-) - \sigma^-(\sigma^-)\T.
\]

\section{Derivations and numerical details for the estimation principles}\label{app:estimation-rev}

\subsection{Dually flat estimation: the EKO principle}

\subsubsection{Shannon running example}\label{app:eko-shannon-example}
For $n = 3$, $t(i) = i$ and the sample mean $\bar x = (0.5,0.2,0.3)$, the empirical moment is $\bar t = 1.8$; solving $\E_{p_\theta}[t] = 1.8$ on the Shannon model gives
\begin{equation}\label{eq:eko-numerisch}
  \hat q_{\mathrm{EKO}} \approx (0.438,\ 0.323,\ 0.238),
\end{equation}
the maximum-entropy distribution on the moment surface.

\subsection{Projectively flat estimation: the WZTY principle}

\subsubsection{Affine reparametrization of the projective geodesics}\label{app:wzty-reparam}
From \eqref{eq:christoffel-proj}, the geodesic equation is
\[
  \ddot q + 2(a \cdot \dot q)\dot q = 0,
  \qquad a = -\lambda\nabla\phi.
\]
Thus the image of a geodesic is a straight line in the $p$-chart, but its affine parameter differs from the Euclidean parameter $\sigma$. Substituting $q(t) = q_0 + \sigma(t)\,u$ reduces the geodesic equation to $\ddot\sigma = 2\lambda\,(\nabla\phi \cdot u)\,\dot\sigma^2$, whose solution is $\dot\sigma \propto e^{2\lambda\phi}$, so
\[
  dt \propto e^{-2\lambda\phi}\,d\sigma
\]
for the mixture connection, and $dt \propto e^{2\lambda\phi}\,d\sigma$ for its dual.

\subsubsection{Derivative of the logarithmic projection score}\label{app:wzty-score}
Let $A(x,q) = 1 + \lambda\nabla\phi(q) \cdot (x - q)$. Gradients are taken in the ambient coordinates of $\R^n$; on the saturated model, stationarity is the vanishing of their component along $V_0$, written below with a multiplier $\nu$ of the constraint $\sum_iq_i = 1$, as in \eqref{eq:lcf-stationarity-fr}. Then
\[
  \nabla_q A = \lambda\big(\nabla^2\phi(q)(x - q) - \nabla\phi(q)\big),
\]
and therefore
\[
  \nabla_qD_\phi^{(\lambda)}(x \mid q)
  = \frac{\nabla^2\phi(q)(x - q) - \nabla\phi(q)}{A(x,q)} + \nabla\phi(q).
\]
Summing over the sample with $w_k = A(x_k,q)^{-1}$, stationarity in the tangent directions reads
\[
  \sum_{k=1}^N w_k(q)\big(\nabla^2\phi(q)(x_k - q) - \nabla\phi(q)\big) + N\,\nabla\phi(q) = \nu\,\one;
\]
division by $W = \sum_kw_k$ gives, with $\hat w_k = w_k/W$,
\[
  \nabla^2\phi(q)\Big(\sum_k\hat w_k(q)\,x_k - q\Big) = \big(1 - N/W(q)\big)\nabla\phi(q) + \frac{\nu}{W}\,\one,
\]
the self-normalized reweighted-mean form used in the text. With $w = A(x,q)^{-1}$, the identity
$1 - w = \lambda\,w\,\nabla\phi(q) \cdot (x - q)$ recombines the single-observation gradient into
\[
  \nabla_qD_\phi^{(\lambda)}(x \mid q) = -\,w\,g(q)\,(x - q),
\]
$g$ the metric \eqref{eq:lambda-metric}; summed over the sample and contracted with $J_\theta\T$, whose columns lie in $V_0$ so that the multiplier drops out, this is \eqref{eq:lambda-score}.

\subsubsection{Interior projection in the running example}\label{app:wzty-interior-example}
At $\lambda = 1$, $n = 3$, $t(i) = i$, the family \eqref{eq:lambda1-family} is $p_\theta = (\tfrac{1}{3} - \theta,\tfrac{1}{3},\tfrac{1}{3} + \theta)$, $|\theta| < 1/3$. For the uniform geometric-mean generator,
\[
  D^{(1)}_{\phi_\pi}(x \mid q)
  = \log\Big(\sum_{i=1}^{3}\pi_i\frac{x_i}{q_i}\Big)
  - \sum_{i=1}^{3}\pi_i\log x_i + \sum_{i=1}^{3}\pi_i\log q_i.
\]
On $x_1 = (0.7,0.1,0.2)$, $x_2 = (0.4,0.3,0.3)$, $x_3 = (0.4,0.2,0.4)$ the scalar objective is strictly convex on the admissible interval and has
\begin{equation}\label{eq:wzty-numeric}
  \hat\theta \approx -0.0755,
  \qquad
  \hat q_{\mathrm{WZTY}} \approx (0.409,\ 0.333,\ 0.258).
\end{equation}
At that point
\[
  w \approx (1.076,0.986,0.958),\qquad
  \sum_k\hat w_kx_k \approx (0.507,0.197,0.296),
\]
and $\E_{\hat q_{\mathrm{WZTY}}}[t] \approx 1.849$, matching neither the raw mean $1.8$ nor the escort mean of the categorical companion.

\subsubsection{Derivative of the categorical escort score}\label{app:wzty-escort-derivative}
For $s_i = \theta \cdot t(i) - \varphi_\lambda(\theta)$,
$\exp_\lambda'(s) = \exp_\lambda(s)^{1-\lambda}$ wherever positive, hence
\[
  \nabla_\theta\log p_\theta(i)
  = p_\theta(i)^{-\lambda}\big(t(i) - \nabla\varphi_\lambda(\theta)\big).
\]
Averaging this score against the empirical frequencies and setting it to zero gives \eqref{eq:begleit-moment}.

\subsubsection{Categorical companion in the running example}\label{app:wzty-categorical-example}
For \eqref{eq:lambda1-family} at $\lambda = 1$ and frequencies $\hat p = (0.5,0.2,0.3)$, \eqref{eq:begleit-moment} becomes
\begin{equation}\label{eq:wzty-numerisch-gl}
  \sum_{i=1}^3\hat p_i\frac{i - 2}{\tfrac{1}{3} + \theta(i - 2)} = 0
  \quad\Longleftrightarrow\quad
  \frac{0.3}{\tfrac{1}{3} + \theta} = \frac{0.5}{\tfrac{1}{3} - \theta}.
\end{equation}
The root is $\hat\theta = -1/12$, hence
\begin{equation}\label{eq:lambda-mle-numerisch}
  \hat q_{\lambda\mathrm{ML}} = (5/12,4/12,3/12) \approx (0.417,0.333,0.250),
\end{equation}
with $\E_{\hat q_{\lambda\mathrm{ML}}}[t] = 11/6$.

\subsection{The Levi-Civita-Fr\'echet principle}

\subsubsection{Proof of the flat-pullback formula}\label{app:lcf-flat-pullback}
Under the hypothesis of \autoref{prop:flat-frechet}, $F$ is a global isometry, so $d(q,x_i) = \|F(q) - F(x_i)\|$. In the variable $\xi = F(q)$ the objective is $\sum_iw_i\|\xi - F(x_i)\|^2$, a strictly convex quadratic minimized at $\xi = \sum_iw_iF(x_i)$. Pulling back through $F^{-1}$ gives \eqref{eq:lcf-coda}.

\subsubsection{Fisher-Rao stationarity and boundary derivative}\label{app:lcf-fr-stationarity}
With \eqref{eq:rao-abstand}, the saturated Fisher-Rao objective is
\[
  F(q) = \frac{4}{N}\sum_{k=1}^N\arccos\Big(\sum_i\sqrt{q_i(x_k)_i}\Big)^2.
\]
Introducing a multiplier $\nu$ for $\sum_iq_i = 1$ and differentiating gives \eqref{eq:lcf-stationarity-fr}. At a face $q_i = 0$ with all observations interior, the factor $\sqrt{(x_k)_i/q_i}$ makes the inward derivative diverge to $-\infty$, so the minimizer cannot lie on the boundary.

\subsubsection{Running Fisher-Rao examples}\label{app:lcf-running-examples}
Solving \eqref{eq:lcf-stationarity-fr} for the interior running sample gives
\begin{equation}\label{eq:lcf-numeric-fr}
  \hat q_{\mathrm{LCF}} \approx (0.503,\ 0.196,\ 0.301),
\end{equation}
whereas the categorical sample with the same arithmetic mean gives
\begin{equation}\label{eq:lcf-numeric-cat}
  \hat q_{\mathrm{LCF}} \approx (0.599,\ 0.135,\ 0.266).
\end{equation}
The two samples have the same arithmetic mean but different Fisher-Rao Fr\'echet means.

\subsection{The symmetrized cost as a loss}

\subsubsection{The barycenter of the running sample}\label{app:jbary-running-example}
For the running sample, direct minimization of \eqref{eq:jbary} gives
\begin{equation}\label{eq:jbary-numeric}
  \hat q_{\mathrm{sym}} \approx (0.504,\ 0.192,\ 0.304),
\end{equation}
at clr-distance approximately $0.006$ from the geometric mean $(0.506,0.191,0.303)$.

\subsection{Conjugacy and the $\lambda$-Gibbs posterior}

\subsubsection{Boundary propriety of the $\lambda$-Gibbs posterior}\label{app:lambda-bayes-propriety}
For the weighted logarithmic generator, near a face $q_i \to 0$ the factor $e^{-N\phi(q)}$ behaves as $q_i^{-N\pi_i}$, while the $N$ data factors in \eqref{eq:lambda-gibbs} contribute $q_i^{N/\lambda}$. The net exponent is $N(1/\lambda - \pi_i)$, so the density vanishes at the face whenever $\lambda\pi_i < 1$. This holds, in particular, for the running uniform-weight examples at $\lambda \le 1$.

\subsubsection{Posterior mode and mean for the running sample}\label{app:lambda-bayes-numerics}
For the running $n = 3$, $N = 3$, $\lambda = 1$ sample the flat-prior posterior \eqref{eq:lambda-gibbs} has the closed rational form
\begin{equation}\label{eq:lambda1-posterior}
  \operatorname{post}(q)\;\propto\;
  \frac{\big(\prod_{i=1}^nq_i\big)^{N(1-1/n)}}{\prod_{k=1}^NP_k(q)},
  \qquad
  P_k(q) = \sum_{i=1}^n x_{k,i}\prod_{j\ne i}q_j,
\end{equation}
a Dirichlet kernel divided by data polynomials of degree $n - 1$. Direct maximization gives the flat-prior mode $(0.517,0.191,0.292)$ and numerical integration the posterior mean $(0.403,0.274,0.323)$; the mode differs from $\bar x = (0.5,0.2,0.3)$ and the posterior mean from both. The randomized-conjugacy identity \eqref{eq:randomisierte-konjugiertheit} is the Feynman parametrization
\begin{equation}\label{eq:feynman}
  \E_{u\sim\Dir(1/\lambda,\dots,1/\lambda)}
  \Big[\Big(\sum_{k=1}^Nu_kB_k\Big)^{-N/\lambda}\Big]
  = \prod_{k=1}^NB_k^{-1/\lambda},
  \qquad B_1,\dots,B_N > 0,
\end{equation}
applied to $B_k = A_k(q)$, an exact identity \citep{feynman1949}.

\subsection{Robustness and efficiency}

\subsubsection{Sensitivity against efficiency along $\beta$}\label{app:rob-eko-numerics}
On the running family at $p = \hat q_{\mathrm{EKO}}$, the gross-error sensitivity of the decoupled density-power estimator falls from $1.885$ at $\beta = 0$, attained at the rarest vertex $e_3$, to $1.778$ at $\beta = 1/2$. The corresponding variance calculation gives $\mathrm{Avar}_0/\mathrm{Avar}_{1/2} \approx 0.95$: increasing $\beta$ from $0$ to $1/2$ reduces the sensitivity, with relative efficiency approximately $0.95$.

\subsubsection{The log-cosh identity}\label{app:rob-wzty-logcosh}
For $n = 2$, uniform weight and $z = \log(q_1/q_2)$, at $q = (1/2,1/2)$ the affine term in the $L$-divergence equals $1$, so
\[
  D_\phi(x \mid q) = -\tfrac{1}{2}\log(4x_1x_2).
\]
Since $4x_1x_2 = \cosh^{-2}(z_x/2)$, this gives \eqref{eq:logcosh} at $z_q = 0$, and perturbation equivariance gives the general difference $z_x - z_q$. Differentiation yields the bounded score $\tfrac{1}{2}\tanh((z_q - z_x)/2)$.

\subsubsection{Influence norms and Hessians for the running samples}\label{app:rob-lcf-numerics}
For the worked interior Fisher-Rao sample, the data lie within $0.42$ of the estimate, the coefficients in \eqref{eq:frechet-hesse} are approximately $(0.986,0.994,0.996)$, and the eigenvalues of $H$ are within two percent of $2$; the influence norms are approximately $(0.41,0.26,0.23)$ and the gross-error sensitivity about $2.25$. For the categorical configuration the coefficients are approximately $(0.838,0.472,0.619)$, the Hessian eigenvalues $(1.657,1.741)$, and the influence norms $(1.58,2.88,2.39)$.

For the categorical sample with equal mass at the three vertices, whose Fr\'echet mean is the barycenter by symmetry, the square-root geometry has distance $2.007$ from the barycenter to a vertex and transverse Jacobi coefficient $0.14$ there, giving $H = 1.14I$ and sensitivity approximately $3.53$; Fisher-Rao on the same sample has distance $1.911$, coefficient $0.68$, $H = 1.68I$ and sensitivity $2.28$. The larger square-root sensitivity comes from both a larger distance and a smaller transverse Hessian eigenvalue.

\subsubsection{The exponential map, Jacobi fields and the distance Hessian}\label{app:rob-geometric-preliminaries}
The exponential map
\begin{equation}\label{eq:exp-abbildung}
  \exp_q: \mathcal V_q \subset T_q\simplex \to \simplex,\qquad v \mapsto \gamma_v(1),
\end{equation}
sends $v$ to the time-one point of the geodesic starting at $q$ with velocity $v$; $\log_q$ is its inverse where the minimizing geodesic is unique. Along a unit-speed geodesic, an infinitesimal normal displacement is a Jacobi field. In dimension two its scalar coefficient solves
\begin{equation}\label{eq:jacobi}
  J''(s) + K(\gamma(s))J(s) = 0,\qquad J(0) = 0,\quad J'(0) = 1,
\end{equation}
and its first positive zero is the first conjugate distance. Positive curvature therefore refocuses nearby geodesics, negative curvature separates them, and the same equation yields the distance Hessian \eqref{eq:hessian-blocks}, whose transverse block is built from the matrix solution of
\begin{equation}\label{eq:matrix-jacobi}
  J''(s) + R(s)\,J(s) = 0,\qquad J(0) = 0,\qquad J'(0) = I,
\end{equation}
along the same geodesic.

\subsubsection{Proof of the existence dichotomy}\label{app:rob-existence-proof}
If one observation lies at infinite distance, every value of the objective contains an infinite summand. Otherwise each extended distance is finite and $1$-Lipschitz on the metric completion; a finite sum of their squares is continuous. On a compact completion it therefore attains its minimum, proving \autoref{prop:existence-dichotomy}.

\subsubsection{Numerical constants of the square-root geometry}\label{app:rob-sqrt-certificates}
For the metric \eqref{eq:wurzel-metrik} with three categories ($n = 3$), the curvature extends to the compactification of \autoref{prop:sqrt-sphere}, with maximum $K_{\max} = 5/4$ at the barycenter and minimum $-1/4$ at the edge midpoints. Along the symmetric rays, joint integration of the geodesic and Jacobi equations gives $dJ'(d)/J(d) = 0.62$ at the edge-midpoint distance $d = 1.251$ and $0.14$ at the vertex distance $d = 2.007$; the latter remains positive beyond the comparison radius $\pi/(2\sqrt{K_{\max}}) = 1.405$. These are numerical values, not certified bounds, and not a global uniqueness theorem for the vertex configuration.

\subsubsection{Symmetrized-loss contamination path}\label{app:rob-sym-numerics}
With a fourth observation $x_4 = C(1,\delta,\delta)$ at weight $1/4$, the clr-norms of the $D^{\mathrm{sym}}$-barycenter and of the geometric mean are
\begin{center}
\begin{tabular}{lrrrrr}
\hline
$\delta$ & $10^{-1}$ & $10^{-2}$ & $10^{-4}$ & $10^{-8}$ & $10^{-12}$\\
\hline
$D^{\mathrm{sym}}$-barycenter & $0.95$ & $1.23$ & $1.30$ & $1.30$ & $1.30$\\
geometric mean & $0.96$ & $1.42$ & $2.35$ & $4.22$ & $6.10$\\
\hline
\end{tabular}
\end{center}
the barycenter approaching a limit near $(0.702,0.122,0.177)$ while the geometric mean diverges logarithmically.

\section{Proofs for the estimation section}\label{app:est-proofs}

Proofs of \autoref{prop:beta-lambda}, \autoref{prop:grad-dist}, \autoref{prop:jbary}, \autoref{prop:bregman-factor}, \autoref{prop:wzty-bounded} and \autoref{prop:sqrt-sphere}, of the Hessian eigenvalues \eqref{eq:sphere-hessian-ev}, and the score of the decoupled density-power loss.

\begin{proof}[Proof of \autoref{prop:beta-lambda}]
By \eqref{eq:lambda-metric}, $\phi_\beta$ is $\lambda$-exponentially concave if and only if
\eqref{eq:beta-lambda-metrik} is positive semidefinite on $T_q\simplex$ (with
$\nabla\phi_\beta = -\xi_\beta(q)$ and
$\nabla^{2}\phi_\beta = -\diag(q_i^{\beta-1})$, which is the metric
claim of the proposition), that is,
\[
  \sum_{i=1}^{n}q_i^{\beta-1}v_i^{2}
  \;\ge\;\lambda\Big(\sum_{i=1}^{n}\xi_\beta(q_i)\,v_i\Big)^{2}
  \quad\text{for all }v \in T_q\simplex.
\]
Since $\sum_iv_i = 0$, the right side is unchanged under
$\xi_\beta(q_i) \mapsto \xi_\beta(q_i) - c$, and Cauchy-Schwarz on the
factorization
$(\xi_\beta(q_i) - c)v_i
= \big(q_i^{(\beta-1)/2}v_i\big)
 \big(q_i^{(1-\beta)/2}(\xi_\beta(q_i) - c)\big)$
bounds it by
$\big(\sum_iq_i^{\beta-1}v_i^{2}\big)
 \sum_iq_i^{1-\beta}\big(\xi_\beta(q_i) - c\big)^{2}$.
Choosing $c$ as the $q^{1-\beta}$-weighted mean of $\xi_\beta(q)$ puts
the equality direction $v_i = q_i^{1-\beta}(\xi_\beta(q_i) - c)$ inside the
tangent space, so the sharp constant at $q$ is exactly $W_\beta(q)$;
taking the supremum over $q$ gives the frontier.

The closed form quoted after the proposition is, with $S = \sum_{i=1}^{n}q_i^{1+\beta}$, escort weights $w_i = q_i^{1+\beta}/S$ and ratios $r_i = p_i/q_i$, on the domain of \eqref{eq:lambda-ldiv},
\begin{equation}\label{eq:beta-lambda-div}
  D^{(\lambda)}_{\phi_\beta}(p \mid q)
  = \frac{1}{\lambda}\log\Big(1 - \frac{\lambda S}{\beta}\big(\E_w[r] - 1\big)\Big)
  + \frac{S}{\beta(1 + \beta)}\big(\E_w[r^{1+\beta}] - 1\big).
\end{equation}
The constant $-1/\beta$ in $\xi_\beta$ contributes nothing, since
$\sum_i(p_i - q_i) = 0$, so
$\nabla\psi_\beta(q) \cdot (p - q)
= \frac{1}{\beta}\big(\sum_iq_i^{\beta}p_i - S\big)
= \frac{S}{\beta}\big(\E_w[r] - 1\big)$,
while
$\psi_\beta(p) - \psi_\beta(q)
= \frac{S}{\beta(1 + \beta)}\big(\E_w[r^{1+\beta}] - 1\big)$,
the terms $-\frac{1}{\beta}\sum_ip_i$ and $-\frac{1}{\beta}\sum_iq_i$ cancelling since both sums equal $1$; inserting both into
\eqref{eq:lambda-ldiv} with $\phi_\beta = -\psi_\beta$ gives
\eqref{eq:beta-lambda-div}.

As $\lambda \to 0$, $\frac{1}{\lambda}\log(1 - \lambda x) \to -x$ turns
\eqref{eq:beta-lambda-div} into
$-\frac{S}{\beta}(\E_w[r] - 1)
+ \frac{S}{\beta(1 + \beta)}(\E_w[r^{1+\beta}] - 1)$, and collecting the
$S$-terms through
$\frac{S}{\beta} - \frac{S}{\beta(1 + \beta)} = \frac{S}{1 + \beta}$ leaves
exactly the summands of the density-power divergence
\begin{equation}\label{eq:dpd}
  D_\beta[p:q] = \sum_{i=1}^n\left\{\frac{1}{1 + \beta}\,q_i^{1+\beta}
  - \frac{1}{\beta}\,p_i\,q_i^{\beta}
  + \frac{1}{\beta(1 + \beta)}\,p_i^{1+\beta}\right\}.
\end{equation}
As $\beta \to 0$,
$\xi_\beta \to \log$, $S \to 1$ and $w \to q$, giving the $\lambda$-tempered Kullback-Leibler divergence
\begin{equation}\label{eq:temperiertes-kl}
  D^{(\lambda)}_{\phi_0}(p \mid q)
  = \frac{1}{\lambda}\log\Big(1 + \lambda\sum_{i=1}^{n}(q_i - p_i)\log q_i\Big)
  + \sum_{i=1}^{n}p_i\log p_i - \sum_{i=1}^{n}q_i\log q_i
\end{equation}
and $W_0(q) = \Var_q(\log q)$; at $\beta = 1$,
$W_1(q) = \min_c\sum_i(q_i - 1 - c)^{2}
= \sum_i\big(q_i - \tfrac{1}{n}\big)^{2} = \sum_iq_i^{2} - \tfrac{1}{n}$, whose
supremum over the simplex, $1 - \tfrac{1}{n}$, is approached at the vertices. Hence at $\beta = 1$ the frontier is $\lambda \le \tfrac{n}{n - 1}$, and at $\beta = 0$ it is $\lambda\,\sup_q\Var_q(\log q) \le 1$. For $n \le 3$ the supremum is below $1$: with $q_1$ the largest coordinate and $r_i = q_i/q_1 \le 1$, choosing $c = \log q_1$ in the variance gives $\Var_q(\log q) \le \sum_iq_i(\log r_i)^2 = q_1\sum_{i \ge 2}r_i(\log r_i)^2$ with $q_1 = 1/(1 + \sum_{i \ge 2}r_i)$, so $\Var_q(\log q) < 1$ as soon as $\sum_{i \ge 2}r_i\big((\log r_i)^2 - 1\big) < 1$; the function $t \mapsto t\big((\log t)^2 - 1\big)$ on $(0,1]$ has maximum $2(1 + \sqrt{2})\,e^{-1-\sqrt{2}} \approx 0.432$ at $t = e^{-1-\sqrt{2}}$, and two such terms stay below $0.864$. For $n \ge 4$ the supremum exceeds $1$: $q = (0.85,0.05,0.05,0.05)$ gives $\Var_q(\log q) \approx 1.023$, and the supremum is nondecreasing in $n$, since letting one coordinate tend to $0$ recovers the value on a face. So at $\lambda = 1$ the Shannon generator is $\lambda$-exponentially concave for $n \le 3$ and not for $n \ge 4$; numerical maximization puts the suprema for $n = 4, 5$ at approximately $1.0235$ and $1.2464$, so for those $n$ it is $\lambda$-exponentially concave up to $\lambda \approx 1/1.0235$ and $1/1.2464$.
\end{proof}

\begin{proof}[Proof of \autoref{prop:grad-dist}]
Vary $q$ along a curve $q_s$ with $\partial_s q_s|_{s=0} = w$; the
corresponding geodesics $\gamma_s$ from $q_s$ to $x$ have energy
$E(\gamma_s) = \int_0^1\|\dot\gamma_s\|_g^2\,dt$, whose first
variation has a boundary term at each endpoint and an interior term that
vanishes identically, because each $\gamma_s$ is itself a geodesic. The
endpoint $x$ is fixed across the family, so only the boundary term at
$q$ survives,
\begin{equation}\label{eq:erste-variation}
  \frac{d}{ds}\bigg|_{s=0} E(\gamma_s)
  = -2\langle w,\dot\gamma(0)\rangle_q.
\end{equation}
Since $E(\gamma) = d(q,x)^2$ for the constant-speed minimizing geodesic
and $\dot\gamma(0) = \log_q(x)$ by definition, this identifies
\eqref{eq:grad-dist-sq}.
\end{proof}

\subsection{Hessian of squared distance under constant curvature}\label{app:sphere-hessian}
\begin{proof}[Proof of \eqref{eq:sphere-hessian-ev}]
Every space of constant curvature $\kappa > 0$ is locally isometric to the
sphere of radius $R = 1/\sqrt{\kappa}$, and the Hessian of
$f_x = d(\cdot,x)^2$ is a local object, so it suffices to compute on that
sphere. Fix $x$ and a point $q$ at distance $d = R\alpha$ from it. Along
the unit-speed geodesic through $q$ and $x$, $f_x(\gamma(s)) = (d + s)^2$,
with second derivative $2$. Along a unit-speed geodesic through $q$
orthogonal to it, the spherical Pythagorean identity gives
$\cos\big(d(\gamma(s),x)/R\big) = \cos\alpha\,\cos(s/R)$, so that
$f_x(\gamma(s)) = R^2\arccos\big(\cos\alpha\cos(s/R)\big)^2$; writing
$u(s)$ for the arccosine, $u(0) = \alpha$, $u'(0) = 0$ and
$u''(0) = \cos\alpha/(R^2\sin\alpha)$, hence
\[
  \frac{d^2}{ds^2}\Big|_{s=0} f_x(\gamma(s))
  = 2R^2\big(u'(0)^2 + u(0)\,u''(0)\big)
  = 2\alpha\cot\alpha
  = 2\sqrt{\kappa}\,d\cot(\sqrt{\kappa}\,d).
\]
By the rotational symmetry of $f_x$ about the radial direction, these
are the eigenvalues of its Hessian: $2$ radially and
$2\sqrt{\kappa}\,d\cot(\sqrt{\kappa}\,d)$ on the orthogonal complement,
both positive while $\sqrt{\kappa}\,d < \pi/2$.
\end{proof}

\begin{proof}[Proof of \autoref{prop:jbary}]
In the chart $\xi = \clr(q)$, with $\xi_k = \clr(x_k)$, the objective is
\[
F(\xi) = \frac{1}{N}\sum_kh(\xi - \xi_k).
\]
Each summand is strictly convex on $V_0$ by \eqref{eq:symcost-hessian},
hence so is $F$. It is also coercive, since \eqref{eq:hilbert-sandwich} bounds
$h(v)$ below by $\big(\max_iv_i - \min_iv_i\big) - 2\log n$, the range being a
norm on $V_0$. A strictly convex coercive function on $V_0$ has exactly one
minimizer, and every point of the chart is interior. Composition with an
injective affine map preserves strict convexity and coercivity.
\end{proof}

\begin{proof}[Proof of \autoref{prop:bregman-factor}]
Averaging the definition of the divergence over $F$ gives
$\E_F D_\psi[X:q] = \E_F D_\psi[X:\mu(F)] + D_\psi[\mu(F):q]$, the
population form of \autoref{lem:two-stage} by the same cancellation; the
first term is free of $q$, so the risk minimizer over any $M$ is
$\Pi_M(\mu(F))$. The mean map is affine in $F$, so along $F_\varepsilon$
its derivative at $\varepsilon = 0$ is $x - \mu(F)$; the chain rule gives
\eqref{eq:ef-affin}.
\end{proof}

\begin{proof}[Proof of \autoref{prop:wzty-bounded}]
Write $S = \sum_i\pi_ir_i$. Since $0 \le \pi_jr_j \le S$, the ratio $r_j/S$
lies in $[0,1/\pi_j]$, so the bracket lies in $[1 - 1/\pi_j,\,1]$; the
endpoints are attained as $x_j \to 0$ (where $r_j/S \to 0$) and $x \to e_j$
(where $S \to \pi_jr_j$).
The tangential score is the projection of the vector $(\psi_j)_j$,
whose norm is bounded by these constants times $\max_j\pi_j/q_j$, and
finiteness of the gross-error sensitivity follows from
\eqref{eq:if-sandwich} with nondegenerate risk Hessian.
\end{proof}

\begin{proof}[Proof of \autoref{prop:sqrt-sphere}]
$S$ is the level set of $s \mapsto \sum_i s_i^4$ at the regular value
$1$, its gradient $4(s_1^3,\dots,s_n^3)$ vanishing only at the origin,
and is closed and bounded. On $S$, Cauchy-Schwarz gives
$\|s\|^2 \le \sqrt{n}\,(\sum_i s_i^4)^{1/2} = \sqrt{n}$, while
$(\sum_i s_i^2)^2 \ge \sum_i s_i^4 = 1$; the conformal factor is smooth and
pinched between $4/\sqrt{n}$ and $4$.

For the isometry, from $p_i = s_i^4$, $dp_i = 4s_i^3\,ds_i$ and
$4\rho\,p_i^{3/2} = 4\|s\|^2 s_i^6$, each summand of
\eqref{eq:wurzel-metrik} transforms as
$dp_i^2/(4\rho\,p_i^{3/2}) = (4/\|s\|^2)\,ds_i^2$; the constraint
$\sum_i p_i = 1$ becomes $\sum_i s_i^4 = 1$, and its differential
$\sum_i dp_i = 0$ becomes $\sum_i s_i^3\,ds_i = 0$, the tangent space of
$S$.

For the geodesics, each sign flip $\sigma_k(s) = (s_1,\dots,-s_k,\dots,s_n)$ preserves
the defining polynomial and $\|s\|$, hence is an isometry of
$(S,\bar g)$, and its fixed-point set, the wall
$\mathcal W_k = \{s \in S: s_k = 0\}$, is totally geodesic. Let $\gamma$ be a
minimizing geodesic between $x,y \in S_+$ (minimizers exist because
the compact $(S,\bar g)$ is complete) and suppose
$\gamma_k(t_0) = 0$ at an interior time. If $\dot\gamma_k(t_0) = 0$, the
geodesic is tangent to $\mathcal W_k$, hence contained in it by
uniqueness of geodesics, contradicting $x \in S_+$. So every zero of
every coordinate is transversal, and there are finitely many. The folded path
$\eta = (|\gamma_1|,\dots,|\gamma_n|)$ lies on $S$
with the same endpoints and, being locally the image of $\gamma$
under compositions of the isometries $\sigma_j$ away from the
crossing times, the same length; it is therefore itself minimizing,
hence a smooth geodesic, yet its $k$-th velocity jumps by
$2|\dot\gamma_k(t_0)| \ne 0$ at $t_0$, a corner. So $\gamma$
never reaches a wall. Interior paths being in particular paths on $S$,
the two distances coincide on $S_+ \times S_+$, realized by the
interior minimizer.

For the completion, the pinching above makes $\bar g$-length and Euclidean length of
curves comparable, so the metric completion of $S_+$ is its Euclidean
closure $\bar S_+$, compact; smoothness of $\bar g$ on
$S \supset \bar S_+$ extends the curvature, and a compact manifold has
positive injectivity radius.
\end{proof}

\subsection{Score of the decoupled density-power loss}\label{app:dpd-score}
For the density-power loss fitted to the undeformed family $E_t$, the per-observation score is
\[
  \psi_\beta(x,\theta) = \zeta_\beta(\theta) - \sum_i x_i\,p_\theta(i)^{\beta}\big(t(i) - \eta(\theta)\big),
  \qquad
  \zeta_\beta(\theta) = \sum_i p_\theta(i)^{1+\beta}\big(t(i) - \eta(\theta)\big),
\]
and differentiating the risk at the model gives the $\beta$-tilted information $J_\beta$ and the influence function \eqref{eq:ef-dpd}.

\section{Minimal-entropy pricing as exponential projection}\label{app:follmer}

A one-period market on a finite sample space $\Omega$ has a reference probability $P$ in the open simplex and $d$ discounted net gains $Y_1,\dots,Y_d$. The set of arbitrage-free pricing measures, $M = \{Q \in \simplex: \E_Q[Y_j] = 0,\ j = 1,\dots,d\}$, is cut out by linear constraints in the mixture coordinates: it is a moment surface $M_\tau$ of \autoref{subsec:eko} at $\tau = 0$, and the fundamental theorem of asset pricing says the market is arbitrage-free if and only if $M$ intersects the open simplex. Minimal-entropy pricing selects from $M$ by relative entropy, $H(Q|P)$ in the notation of F\"ollmer and Schied and $D_{\mathrm{KL}}(Q\|P)$ in ours:
\begin{equation}\label{eq:memm}
\hat Q
= \arg\min_{Q\in M} D_{\mathrm{KL}}(Q\|P),
\qquad
\frac{d\hat Q}{dP}
= \frac{e^{\lambda^{\ast}\cdot Y}}{\E_{P}[e^{\lambda^{\ast}\cdot Y}]},
\qquad
\E_{\hat Q}[Y] = 0.
\end{equation}
Corollary~3.25 of F\"ollmer and Schied \citep{follmerschied2016} gives existence, uniqueness and the exponential form of the minimizer, under the name entropy-minimizing risk-neutral measure. In the terms of \autoref{subsec:eko}, $\hat Q$ is the e-projection of $P$ onto the m-flat surface $M$; it belongs to the exponential family generated by $Y$ with reference measure $P$; the multiplier $\lambda^\ast$ is fixed by the moment condition $\E_{\hat Q}[Y] = 0$; and the Pythagorean relation of \autoref{thm:pyth}, $D_{\mathrm{KL}}(Q\|P) = D_{\mathrm{KL}}(Q\|\hat Q) + D_{\mathrm{KL}}(\hat Q\|P)$ for every $Q \in M$, is display (3.12) in the proof of their Theorem~3.24. The minimal entropy also has a Legendre-dual representation: since $\E_{\hat Q}[Y] = 0$, it is minus the minimal cumulant, $D_{\mathrm{KL}}(\hat Q\|P) = -\min_\lambda\log\E_P[e^{\lambda\cdot Y}]$, the conjugate of the log-partition function at $\tau = 0$. The selected pricing measure depends on the divergence: relative entropy gives the e-projection, with exponential densities and moment matching, and another divergence of \autoref{sec:Estimators} can select a different measure from the same $M$.

\section{Numerical details for the applications}\label{app:applications}

\subsection{The curvature exhibit}\label{app:cid-curvature}
Through the isometry \eqref{eq:sphere-iso}, $\xi = 2\sqrt{\omega}$ on the sphere of radius $2$,
\begin{equation}\label{eq:fr-transport}
\begin{gathered}
\cos\theta = \tfrac{1}{4}\langle \xi_x,\xi_y\rangle,\qquad
\log_{x}(y) = \frac{\theta}{\sin\theta}\left(\xi_y - \cos\theta\,\xi_x\right),\\
\tau_{x\to y}(w) = w - \frac{\langle \xi_y,w\rangle}{4 + \langle \xi_x,\xi_y\rangle}\left(\xi_x + \xi_y\right),
\end{gathered}
\end{equation}
with $\exp_x(v) = \cos(|v|/2)\,\xi_x + 2\sin(|v|/2)\,v/|v|$ and shares recovered as $\omega_i = \xi_i^2/4$. For the projection transport of Zhou et al.~\citep{zhou2025}, write the control direction at $x = \sqrt{\omega_{00}}$ as $a\,u + w$, with $u$ the unit tangent toward the treated baseline, $w$ orthogonal to the great circle through the two base points and $\varphi$ the angle between them; projection returns $a\cos\varphi\,u' + w$, rescaled to the length of $a\,u + w$, where Levi-Civita transport returns $a\,u' + w$, and the rescaled vector equals $a\,u' + w$ only if $a = 0$, $w = 0$ or $\varphi = 0$.

Control cells as in \autoref{tab:cid-parallelisms}; treated cells $\omega_{10} = (0.20,0.60,0.20)$, $\omega^{I}_{11} = (0.15,0.55,0.30)$, chosen so that the three base cells span a proper geodesic triangle while every transport stays interior. The square-root transport is computed by integrating the geodesic and parallel-transport equations in the fourth-root chart; the loop rotation reproduces the Gauss-Bonnet angle-sum excess to $4 \times 10^{-13}$ for a curvature field with range $[-\tfrac{1}{4},\tfrac{5}{4}]$. In shares, the Fisher-Rao direct counterfactual is $(0.0058,0.5323,0.4619)$ against the detour $(0.0040,0.5603,0.4357)$.

\begin{table}[h]
\centering
\begin{tabular}{p{0.72\textwidth}r}
\hline
quantity & value\\
\hline
Fisher-Rao loop rotation $=$ spherical excess & $4.046^\circ$\\
square-root-metric loop rotation $=$ Gauss-Bonnet excess & $4.273^\circ$\\
path dependence of the counterfactual, Fisher-Rao LC (FR distance, direct vs.\ detour) & $0.0606$\\
path dependence of the counterfactual, square-root LC (FR distance) & $0.0659$\\
direction gap at $\omega_{10}$, projection vs.\ Levi-Civita & $4.163^\circ$\\
for comparison: FR distance between the e- and the LC-counterfactual & $0.2892$\\
\hline
\end{tabular}
\caption{The curvature exhibit on the cells of \autoref{app:cid-curvature}.}
\label{tab:cid-curvature}
\end{table}

For the treated baselines $\omega_{10} = (0.15,0.55,0.30)$ and $(0.10,0.60,0.30)$, both interior, the Fisher-Rao exponential map crosses the wall $\omega_1 = 0$; at the second design the projection transport crosses as well, while the square-root transport stays interior with chart coordinates bounded away from the wall.

\bibliographystyle{plainnat}
\bibliography{references}

\end{document}